\documentclass[a4paper,11pt]{article}

\usepackage{fullpage,amsthm,enumitem,mathtools,amssymb}
\newtheoremstyle{break}
  {0em}{0em}%
  {\itshape}{}%
  {\bfseries}{}%
  {\newline}{}%
\theoremstyle{break}
\usepackage[linewidth=0.5pt,backgroundcolor=frame,linecolor=white]{mdframed}
\newmdtheoremenv[linewidth=0.5pt,backgroundcolor=frame,linecolor=white,innertopmargin=0.4em,innerbottommargin=0.4em,innerleftmargin=0.4em,innerrightmargin=0.4em]{definition}{Definition}
\newmdtheoremenv[linewidth=0.5pt,backgroundcolor=frame,linecolor=white,innertopmargin=0.4em,innerbottommargin=0.4em,innerleftmargin=0.4em,innerrightmargin=0.4em]{algorithm}{Algorithm}
\definecolor{frame}{rgb}{0.9,0.9,0.9}
\definecolor{colorgray}{rgb}{0.8,0.8,0.8}
\usepackage{enumerate}
\usepackage{amsmath,subfig,float,graphicx,multirow}
\usepackage{hyperref}
\hypersetup{colorlinks = true, linktocpage = true, bookmarksopen = true, linkcolor = blue, urlcolor=blue, citecolor = blue, allcolors=blue}
\graphicspath{{"./Figures/"}}
\usepackage[numbers,sort&compress]{natbib}
\usepackage{amsthm}
\theoremstyle{definition}
\newtheorem{proposition}{Proposition}
\renewenvironment{proof}{\noindent{\bfseries Proof.}}

\makeatletter
\newcommand{\redub}{}
\def\redub#1{%
  \@ifnextchar_%
    {\@redub{#1}}
    {\@latex@warning{Missing argument for \string\redub}\@redub{#1}_{}}%
}
\def\@redub#1_#2{%
    \underbracket[0.5pt]{\color{black}#1}_{\displaystyle\color{black} #2}%
}
\makeatother

\begin{document}

\title{Semi-analytical solution of multilayer diffusion problems with time-varying boundary conditions and general interface conditions}
\date{}
\author{Elliot J. Carr\footnote{Corresponding author: \href{mailto:elliot.carr@qut.edu.au}{elliot.carr@qut.edu.au}.} {} and Nathan G. March\\ School of Mathematical Sciences, Queensland University of Technology (QUT),\\ Brisbane, Australia.}
\maketitle

\begin{abstract}
We develop a new semi-analytical method for solving multilayer diffusion problems with time-varying external boundary conditions and general internal boundary conditions at the interfaces between adjacent layers. The convergence rate of the semi-analytical method, relative to the number of eigenvalues, is investigated and the effect of varying the interface conditions on the solution behaviour is explored. Numerical experiments demonstrate that solutions can be computed using the new semi-analytical method that are more accurate and more efficient than the unified transform method of Sheils [\textit{Appl.~Math.~Model.}, 46:450--464, 2017]. Furthermore, unlike classical analytical solutions and the unified transform method, only the new semi-analytical method is able to correctly treat problems with both time-varying external boundary conditions and a large number of layers. The paper is concluded by replicating solutions to several important industrial, environmental and biological applications previously reported in the literature, demonstrating the wide applicability of the work.\\
\\
\textbf{Keywords:}~
~multilayer diffusion; semi-analytical solution; transient boundary
conditions; general interface conditions; partition coefficient; jump
conditions
\end{abstract}

\section{Introduction}

Mathematical models of diffusion in layered materials arise in many industrial, environmental, biological and medical applications, such as heat conduction in composite materials \cite{monte_2000,mikhailov_1983,mulholland_1972}, transport of contaminants, chemicals and gases in layered porous media \cite{liu_1998,yates_2000}, brain tumour growth \cite{asvestas_2014,mantzavinos_2014}, heat conduction through skin \cite{simpson_2017}, transdermal drug delivery \cite{todo_2013,pontrelli_2007} and greenhouse gas emissions \cite{liu_2008}. Another important application is in the field of multiscale modelling: for a large number of layers, layered diffusion is one of the simplest examples of a multiscale problem and is ideal for prototyping macroscopic and multiscale modelling approaches \cite{carr_2016a,carr_2016b}. Analytical solutions to such problems are highly valuable as they provide a greater level of insight into the solution behaviour and can be used to benchmark numerical solutions.

The most popular analytical solution approach for multilayer diffusion is classical separation of variables (see, e.g., \cite{mikhailov_1983,hickson_2009a,trefry_1999}).  By assuming a separated solution in each layer, one immediately finds that the eigenfunctions  are coupled via the internal boundary conditions (BCs) at the interfaces between adjacent layers. Together with the external BCs, a system of algebraic equations is obtained, linear in the (unknown) eigenfunction coefficients. With the requirement that this linear system possess non-trivial solutions, one obtains a transcendental equation satisfied by the unknown eigenvalues formulated by setting the determinant of the coefficient matrix of the linear system equal to zero. Since computing the matrix determinant is numerically unstable for large matrices, using this method to compute the analytical solution performs poorly for a large number of layers as either erroneous eigenvalues/roots are returned during the solution procedure or eigenvalues/roots are skipped altogether (as reported by \citet{carr_2016a}).

To overcome these issues, \citet{carr_2016a} recently developed a semi-analytical solution approach for multilayer diffusion based on the Laplace transform and an appropriately-defined orthogonal eigenfunction expansion. The attractiveness of this approach is that the solution formulas involve a local set of eigenvalues in each layer satisfying simple transcendental equations, resembling those of the single layer problem, which in most cases can be solved explicitly. As a result, the solution performs well for a large number of layers \cite{carr_2016a} with the approach classified as \textit{semi-analytical} since computing the inverse Laplace transforms appearing in the solution formulas are carried out numerically.

Recently, \citet{sheils_2016} applied the unified transform method, initially proposed by \citet{fokas_1997}, to the layered diffusion problem and compared the approach to the semi-analytical method of \citet{carr_2016a}. Sheils reported that her approach is more accurate, but less efficient and also faulty near the boundaries whenever nonhomogeneous external BCs were applied. We note that her assessment of the accuracy was based on using the default value of 50 terms/eigenvalues in the solution expansions in \citet{carr_2016a}'s code. \citet{sheils_2016} also remarked that a notable difference between the two methods is that only the unified transform method is applicable to time-dependent external BCs.

In this paper, several key contributions to the literature on multilayer diffusion are presented. Specifically, we:
\begin{itemize}
\item[(i)] develop a new semi-analytical solution approach for time-dependent external BCs and a general set of internal BCs (Section \ref{sec:semi-analytical});
\item[(ii)] study the convergence rate of the new semi-analytical method in (i) and compare it to the classical analytical method (Section \ref{sec:convergence});
\item[(iii)] carry out a more comprehensive comparison between our semi-analytical method and Sheils' \cite{sheils_2016} unified transform method, in terms of accuracy and efficiency, by performing an investigation into the effect of the number of terms/eigenvalues on the solution accuracy (Section \ref{sec:compare_Sheils});
\item[(iii)] explore the effect that changing the interface conditions has on the solution behaviour (Section \ref{sec:interface_conditions});
\item[(iv)] extend the classical analytical solution approach to a general set of internal BCs by proposing the correct form of the weight function in the orthogonality condition (Appendix \ref{app:analytical_solution}).
\end{itemize}
The treatment of time-dependent external BCs addresses a deficiency of \citet{carr_2016a}'s semi-analytical method, as reported by \citet{sheils_2016}, and allows for application of the method to a wider range of problems, for example, contaminant transport modelling in layered porous media involving a time-varying inlet concentration \cite{liu_1998}. Our method is also not faulty at the end points as is the case for Sheil's \cite{sheils_2016} unified transform method whenever nonhomogeneous external BCs are applied. Moreover, general interface conditions permit additional features to be incorporated in the layered diffusion model, such as the volumetric heat capacity in heat conduction problems \cite{hickson_2011} or the sorption coefficient and partitioning phenomena in chemical transport problems \cite{trefry_1999}.

The remaining sections of this paper are organised in the following manner. Section \ref{sec:multilayer_problem} formulates the multilayer diffusion problem considered in this work while the derivation and implementation of the semi-analytical solution method is outlined in Section \ref{sec:semi-analytical}. Numerical experiments are reported in Section \ref{sec:results}. The paper then concludes with a summary and overview of key findings.

\section{Multilayer diffusion problem}
\label{sec:multilayer_problem}
\begin{figure}[h]
\centering
\includegraphics[width=0.9\textwidth]{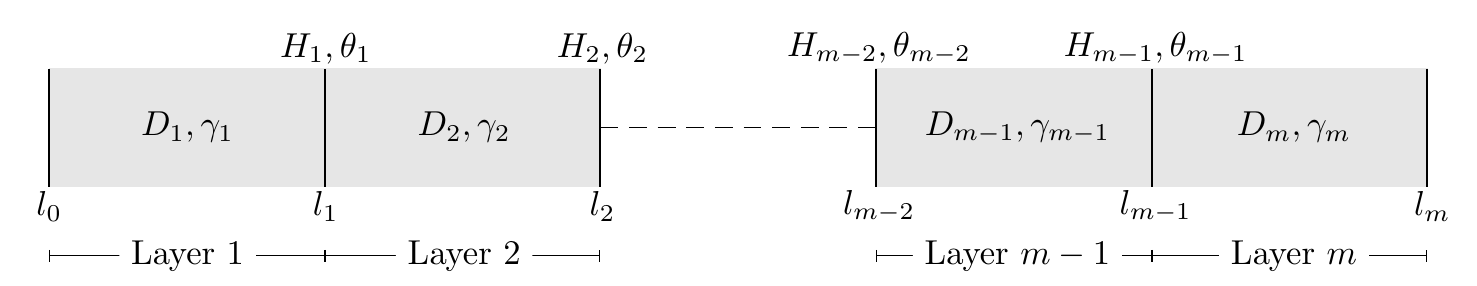}
\caption{One-dimensional layered medium. The coefficients $D_{i}$ and $\gamma_{i}$ are constant in each layer ($i=1,\hdots,m$) while the contact transfer coefficient $H_{i}$ and partition coefficient $\theta_{i}$ are constants defined at the interfaces ($i=1,\hdots,m-1$).}
\label{fig:layered_medium}
\end{figure}

The one-dimensional multilayer diffusion problem is formulated as follows. Consider a diffusion process on the interval $[l_{0},l_{m}]$, which is partitioned into $m$ subintervals (see Figure \ref{fig:layered_medium}):
\begin{align}
l_{0} < l_{1} < \hdots < l_{m-1} < l_{m}.
\end{align}
{The resulting domain is represented using the notation $[l_0,l_1,\hdots,l_{m-1},l_m]$.} On each subinterval (layer) $(l_{i-1},l_{i})$, that is for all $i = 1,\hdots,m$, we define the diffusion equation:
\begin{gather}
\label{eq:original_equation}
\frac{\partial u_{i}}{\partial t} = D_{i}\frac{\partial^{2}u_{i}}{\partial x^{2}},
\end{gather}
where $x\in (l_{i-1},l_{i})$, $u_{i}(x,t)$ is the solution (temperature, concentration, etc.) at position $x$ and time $t$ in the $i$th layer, and $D_{i} > 0$ is the diffusion coefficient in the $i$th layer. The initial conditions are given by 
\begin{gather}
\label{eq:original_ic}
u_{i}(x,0) = f_{i}(x),
\end{gather}
for all $i =1,\hdots,m$, and the external BCs are defined as
\begin{gather}
\label{eq:original_bc1}
a_{L}u_{1}(l_{0},t) - b_{L}\frac{\partial u_{1}}{\partial x}(l_{0},t) = g_{0}(t),\\
\label{eq:original_bc2}
a_{R}u_{m}(l_{m},t) + b_{R}\frac{\partial u_{m}}{\partial x}(l_{m},t) = g_{m}(t),
\end{gather}
for $t > 0$, where the coefficients $a_{L}$, $b_{L}$, $a_{R}$ and $b_{R}$ are non-negative constants satisfying $a_{L} + b_{L} > 0$ and $a_{R} + b_{R} > 0$, and $g_{0}(t)$ and $g_{m}(t)$ are specified time-dependent functions. Note that the restriction on the coefficients is a classical constraint placed on the single-layer problem ($m = 1$) \cite{trim_1990}, which ensures all of the eigenvalues (see Appendix \ref{app:analytical_solution}) are non-negative and the solution remains bounded as $t\rightarrow\infty$.

Closing the problem requires equations (\ref{eq:original_equation})--(\ref{eq:original_bc2}) to be coupled with appropriate internal BCs at the interfaces between adjacent layers: $x = l_{i}$ for all $i = 1,\hdots,m-1$. Typically, continuity of the solution and the diffusive flux is implicitly assumed at each interface, that is,
\begin{subequations}
\label{eq:ic_implicit}
\begin{gather}
\label{eq:type1_1}
u_{i}(l_{i},t) = u_{i+1}(l_{i},t),\\
\label{eq:type1_2}
D_{i}\frac{\partial u_{i}}{\partial x}(l_{i},t) = D_{i+1}\frac{\partial u_{i+1}}{\partial x}(l_{i},t),
\end{gather}
\end{subequations}
for $t > 0$ and for all $i = 1,\hdots,m-1$. In this work, in addition to the above interface conditions (\ref{eq:ic_implicit}), we also consider the following choices of interface conditions.

\textbf{Perfect contact conditions.} Consider the following interface condition at the $i$th interface: 
\begin{subequations}
\label{eq:ic_perfect}
\begin{gather}
\label{eq:ic_perfect_1}
u_{i}(l_{i},t) = u_{i+1}(l_{i},t),\\
\label{eq:ic_perfect_2}
\gamma_{i}\frac{\partial u_{i}}{\partial x}(l_{i},t) = \gamma_{i+1}\frac{\partial u_{i+1}}{\partial x}(l_{i},t),
\end{gather}
\end{subequations}
for $t > 0$, where $\gamma_{i} > 0$. These interface conditions generalise (\ref{eq:ic_implicit}) and permit a wider array of problems to be considered as $\gamma_{i}$ and/or $\gamma_{i+1}$ can be different from $D_{i}$ and $D_{i+1}$, respectively. For example, in heat transfer, where $D_{i}$ is the thermal diffusivity $D_{i} := k_{i}/(\rho_{i}{c_{p}}_{i})$ ($k_{i}$ and $\rho_{i} {c_{p}}_{i}$ are the thermal conductivity and volumetric heat capacity in layer $i$, respectively), the interface conditions (\ref{eq:ic_perfect}) allow continuity of the heat flux to be imposed at the interface by setting $\gamma_{i} = k_{i}$.

\textbf{Jump conditions.} Consider the following interface condition at the $i$th interface: 
\begin{subequations}
\label{eq:ic_jump}
\begin{gather}
\label{eq:ic_jump_1}
\gamma_{i}\frac{\partial u_{i}}{\partial x}(l_{i},t) = H_{i}(u_{i+1}(l_{i},t) - u_{i}(l_{i},t)),\\
\label{eq:ic_jump_2}
\gamma_{i}\frac{\partial u_{i}}{\partial x}(l_{i},t) = \gamma_{i+1}\frac{\partial u_{i+1}}{\partial x}(l_{i},t),
\end{gather}
\end{subequations}
for $t > 0$, where $H_{i}>0$ is the contact transfer coefficient at $x = l_{i}$. These interfaces conditions are a more general form of the perfect contact conditions (\ref{eq:ic_perfect}) since dividing equation (\ref{eq:ic_jump_1}) by $H_{i}$ yields (\ref{eq:ic_perfect_1}) in the limit $H_{i}\rightarrow\infty$. For a finite value of $H_{i}$, the interface conditions (\ref{eq:ic_jump}) produce a discontinuity or jump in the solution field at the interfaces. This is useful in applications involving contact resistance at the interfaces, for example in heat transfer, where a very thin resistive (low conductive) layer causes a jump in the temperature \cite{hickson_2011}.

\textbf{Partition conditions.} Consider the following interface condition at the $i$th interface:
\begin{subequations} 
\label{eq:ic_partition}
\begin{gather}
\label{eq:ic_partition_1}
u_{i}(l_{i},t) = \theta_{i}u_{i+1}(l_{i},t),\\
\label{eq:ic_partition_2}
\gamma_{i}\frac{\partial u_{i}}{\partial x}(l_{i},t) = \gamma_{i+1}\frac{\partial u_{i+1}}{\partial x}(l_{i},t),
\end{gather}
\end{subequations}
for $t > 0$, where $\theta_{i}>0$ is the \textit{partition coefficient} at $x = l_{i}$. Again, these interface conditions are a more general form of the perfect contact conditions (\ref{eq:ic_perfect}), which are given by the special case: $\theta_{i} = 1$. The interface conditions (\ref{eq:ic_partition}) maintain a constant ratio between the discontinuous solution values at the interface. This phenomena is common in applications such as analyte transport in porous media \cite{trefry_1999} and drug release from multi-layer capsules \cite{kaoui_2017}.

\textbf{General interface conditions.} Each of the above sets of internal BCs can be expressed in the general form:
\begin{subequations}
\label{eq:ic_general}
\begin{gather}
\label{eq:ic_general1}
\gamma_{i}\frac{\partial u_{i}}{\partial x}(l_{i},t) = H_{i}(\theta_{i} u_{i+1}(l_{i},t) - u_{i}(l_{i},t)),\\
\label{eq:ic_general2}
\gamma_{i}\frac{\partial u_{i}}{\partial x}(l_{i},t) = \gamma_{i+1}\frac{\partial u_{i+1}}{\partial x}(l_{i},t),
\end{gather}
\end{subequations}
as setting $\gamma_{i} = D_{i}$, $\gamma_{i+1} = D_{i+1}$ and $\theta_{i} = 1$, and taking $H_{i}\rightarrow\infty$ produces (\ref{eq:ic_implicit}); setting $\theta_{i} = 1$ and taking $H_{i}\rightarrow\infty$ yields the perfect contact conditions (\ref{eq:ic_perfect}); setting $\theta_{i} = 1$ produces the jump conditions (\ref{eq:ic_jump}); and taking $H_{i}\rightarrow\infty$ yields the partition conditions (\ref{eq:ic_partition}). 

\section{Semi-analytical solution}
\label{sec:semi-analytical}
We now develop our new semi-analytical method for solving the multilayer diffusion problem with time-varying external BCs (\ref{eq:original_equation})--(\ref{eq:original_bc2}) subject to the general internal BCs (\ref{eq:ic_general}).

\subsection{Reformulation of problem}
\label{sec:reformulation}
Define the interface functions \cite{carr_2016a}:
\begin{align}
\label{eq:diffusion_problem_transient_bc_gi}
g_{i}(t) := \gamma_{i}\frac{\partial u_{i}}{\partial x}(l_{i},t) = \gamma_{i+1}\frac{\partial u_{i+1}}{\partial x}(l_{i},t),
\end{align}
for all $i = 1,\hdots,m-1$. Equations (\ref{eq:original_equation})--(\ref{eq:original_bc2}) and (\ref{eq:ic_general2}) can now be reformulated as a sequence of single layer problems \cite{rodrigo_2016}:
\begin{itemize}[leftmargin=1em]
\item First Layer ($i = 1$):
\begin{subequations}
\label{eq:u_first_layer}
\begin{gather}
\frac{\partial u_{1}}{\partial t} = D_{1}\frac{\partial^{2}u_{1}}{\partial x^{2}},\\
\label{eq:u_first_layer_bc1}
u_{1}(x,0) = f_{1}(x),\quad
a_{L}u_{1}(l_{0},t) - b_{L}\frac{\partial u_{1}}{\partial x}(l_{0},t) = g_{0}(t),\quad
\gamma_{1}\frac{\partial u_{1}}{\partial x}(l_{1},t) = g_{1}(t).\end{gather}
\end{subequations}
\item Middle Layers ($i = 2,\hdots,m-1$):
\begin{subequations}
\label{eq:u_middle_layers}
\begin{gather}
\frac{\partial u_{i}}{\partial t} = D_{i}\frac{\partial^{2}u_{i}}{\partial x^{2}},\\
\label{eq:u_middle_layers_bc1}
u_{i}(x,0) = f_{i}(x),\quad
\gamma_{i}\frac{\partial u_{i}}{\partial x}(l_{i-1},t) = g_{i-1}(t),\quad
\gamma_{i}\frac{\partial u_{i}}{\partial x}(l_{i},t) = g_{i}(t).
\end{gather}
\end{subequations}
\item Last Layer ($i = m$):
\begin{subequations}
\label{eq:u_last_layer}
\begin{gather}
\frac{\partial u_{m}}{\partial t} = D_{m}\frac{\partial^{2}u_{m}}{\partial x^{2}},\\
\label{eq:u_last_layer_bc1}
\hspace*{-0.5cm}u_{m}(x,0) = f_{m}(x),\quad
\gamma_{m}\frac{\partial u_{m}}{\partial x}(l_{m-1},t) = g_{m-1}(t),\quad
a_{R}u_{m}(l_{m},t) + b_{R}\frac{\partial u_{m}}{\partial x}(l_{m},t) = g_{m}(t).\hspace{-0.1cm}
\end{gather}
\end{subequations}
\end{itemize}
Clearly, the solution of each of the above problems will involve the interface functions (\ref{eq:diffusion_problem_transient_bc_gi}), which are unknown. The general idea is therefore to solve the above single layer problems subject to the constraint that the solutions satisfy the, as yet unused, interface condition (\ref{eq:ic_general1}) \cite{carr_2016a,rodrigo_2016}, which can be rewritten in terms of $g_{i}(t)$ as follows
\begin{align}
\label{eq:original_constraint}
\frac{1}{H_{i}}g_{i}(t) = \theta_{i}u_{i+1}(l_{i},t) - u_{i}(l_{i},t),
\end{align}
for all $i = 1,\hdots,m-1$. In summary, solving the single layer problems (\ref{eq:u_first_layer})--(\ref{eq:u_last_layer}) subject to the constraint (\ref{eq:original_constraint}) is equivalent to solving the multilayer diffusion problem described by equations (\ref{eq:original_equation})--(\ref{eq:original_bc2}) and (\ref{eq:ic_general}).

\subsection{Solution of the single layer problems}
To solve each of the single-layer problems (\ref{eq:u_first_layer})--(\ref{eq:u_last_layer}), we introduce the substitution: 
\begin{align}
\label{eq:substitution_u}
u_{i}(x,t) = w_{i}(x,t) + v_{i}(x,t),
\end{align}
for all $i = 1,\hdots,m$, where $w_{i}(x,t)$ is chosen so that $v_{i}(x,t)$ satisfies homogeneous versions of the BCs given in equations (\ref{eq:u_first_layer_bc1}), (\ref{eq:u_middle_layers_bc1}) and (\ref{eq:u_last_layer_bc1}). For example, in the first layer, $w_{1}(x,t)$ satisfies the BCs:
\begin{gather}
\label{eq:w1_bc1}
a_{L}w_{1}(l_{0},t) - b_{L}\frac{\partial w_{1}}{\partial x}(l_{0},t) = g_{0}(t),\quad
\gamma_{1}\frac{\partial w_{1}}{\partial x}(l_{1},t) = g_{1}(t).\end{gather}
If $a_{L}\neq 0$, then one choice\footnote{{This choice, while not unique, does not effect the uniqueness of the final solution as the initial conditions and source/sink terms are appropriately corrected in (\ref{eq:fi_modified})--(\ref{eq:Gi_source}).}} for $w_{1}$ is a linear function in $x$ with time-dependent coefficients. Substituting $w_{1}(x,t) = A(t) + B(t)x$ into (\ref{eq:w1_bc1}) and solving for $A(t)$ and $B(t)$ gives
\begin{align*}
A(t) = \frac{g_{0}(t)}{a_{L}} + \left(\frac{b_{L}}{a_{L}}-l_{0}\right)\frac{g_{1}(t)}{\gamma_{1}},\quad B(t) = \frac{g_{1}(t)}{\gamma_{1}},
\end{align*}
and hence
\begin{align}
\label{eq:w1}
w_{1}(x,t) = \underbracket[0.5pt]{\frac{1}{a_{L}}}_{\psi_{1,1}(x)}\hspace*{-0.2cm}g_{0}(t) + \underbracket[0.5pt]{\frac{a_{L}(x-l_{0})+b_{L}}{\gamma_{1}a_{L}}}_{\psi_{1,2}(x)}g_{1}(t).
\end{align}
Repeating this process, one sees that, in general, the function $w_{i}(x,t)$ can be expressed as a linear combination of the unknown interface functions at $x = l_{i-1}$ and $x = l_{i}$, that is
\begin{align}
\label{eq:wi}
w_{i}(x,t) = g_{i-1}(t)\psi_{i,1}(x) + g_{i}(t)\psi_{i,2}(x),
\end{align}
for all $i = 1,\hdots,m$. The functions $\psi_{i,1}(x)$ and $\psi_{i,2}(x)$ are identified in equation (\ref{eq:w1}) for $i = 1$ and $a_{L}\neq 0$. The remaining cases are summarised in Appendix \ref{app:psi_functions}. Substituting (\ref{eq:substitution_u}) into (\ref{eq:u_first_layer})--(\ref{eq:u_last_layer}), we see that the functions $v_{i}(x,t)$ satisfy the following problems with homogeneous BCs:
\begin{itemize}[leftmargin=1em]
\item First Layer ($i = 1$):
\begin{subequations}
\label{eq:v_first_layer}
\begin{gather}
\frac{\partial v_{1}}{\partial t} = D_{1}\frac{\partial^{2}v_{1}}{\partial x^{2}} + G_{1}(x,t),\\
v_{1}(x,0) = \tilde{f}_{1}(x),\quad
a_{L}v_{1}(l_{0},t) - b_{L}\frac{\partial v_{1}}{\partial x}(l_{0},t) = 0,\quad
\gamma_{1}\frac{\partial v_{1}}{\partial x}(l_{1},t) = 0.
\end{gather}
\end{subequations}
\item Middle Layers ($i = 2,\hdots,m-1$):
\begin{subequations}
\label{eq:v_middle_layers}
\begin{gather}
\label{eq:pde_vi}
\frac{\partial v_{i}}{\partial t} = D_{i}\frac{\partial^{2}v_{i}}{\partial x^{2}} + G_{i}(x,t),\\
\label{eq:vi_ic}
v_{i}(x,0) = \tilde{f}_{i}(x),\quad
\gamma_{i}\frac{\partial v_{i}}{\partial x}(l_{i-1},t) = 0,\quad
\gamma_{i}\frac{\partial v_{i}}{\partial x}(l_{i},t) = 0.
\end{gather}
\end{subequations}
\item Last Layer ($i = m$):
\begin{subequations}
\label{eq:v_last_layer}
\begin{gather}
\frac{\partial v_{m}}{\partial t} = D_{m}\frac{\partial^{2}v_{m}}{\partial x^{2}} + G_{m}(x,t),\\
v_{m}(x,0) = \tilde{f}_{m}(x),\quad
\gamma_{m}\frac{\partial v_{m}}{\partial x}(l_{m-1},t) = 0,\quad
a_{R}v_{m}(l_{m},t) + b_{R}\frac{\partial v_{m}}{\partial x}(l_{m},t) = 0.
\end{gather}
\end{subequations}
\end{itemize}
The modified initial conditions and the source/sink terms are defined as:
\begin{gather}
\label{eq:fi_modified}
\tilde{f}_{i}(x) = f_{i}(x) - \left[g_{i-1}(0)\psi_{i,1}(x) + g_{i}(0)\psi_{i,2}(x)\right],\\
\label{eq:Gi_source}
G_{i}(x,t) = D_{i}\left[g_{i-1}(t)\psi_{i,1}''(x) + g_{i}(t)\psi_{i,2}''(x)\right] - \left[g_{i-1}'(t)\psi_{i,1}(x) + g_{i}'(t)\psi_{i,2}(x)\right],
\end{gather}
for all $i = 1,\hdots,m$. The solution of each of the single layer problems (\ref{eq:v_first_layer})--(\ref{eq:v_last_layer}) can be expressed as an eigenfunction expansion \cite{trim_1990}:
\begin{align*}
v_{i}(x,t) = \sum_{n=0}^{\infty}c_{i,n}(t)\widehat{\phi}_{i,n}(x),
\end{align*}
with coefficients and orthonormal eigenfunctions defined as
\begin{gather*}
c_{i,n}(t) = \int_{l_{i-1}}^{l_{i}} v_{i}\widehat{\phi}_{i,n}\,dx,\quad\widehat{\phi}_{i,n} = \frac{\phi_{i,n}}{\|\phi_{i,n}\|_{2}},\quad \|\phi_{i,n}\|_{2}^{2} = \int_{l_{i-1}}^{l_{i}} \phi_{i,n}^{2}\,dx.
\end{gather*}
The eigenvalues ($\lambda_{i}$) and non-normalized eigenfunctions ($\phi_{i}$) satisfy:
\begin{gather*}
-\phi_{i}'' = \lambda_{i}^{2}\phi_{i},
\end{gather*}
for $x\in(l_{i-1},l_{i})$ and $i=1,\hdots,m$, subject to the following BCs:
\begin{itemize}
\item First Layer ($i=1$): $a_{L}\phi_{1}(l_{0}) - b_{L}\phi_{1}'(l_{0}) = 0$ and $\phi_{1}'(l_{1}) = 0$.
\item Middle Layers ($i=2,\hdots,m-1$): $\phi_{i}'(l_{i-1}) = 0$ and $\phi_{i}'(l_{i}) = 0$.
\item Last Layer ($i=m$): $\phi_{m}'(l_{m-1}) = 0$ and $a_{R}\phi_{m}(l_{m})+b_{R}\phi_{m}'(l_{m}) = 0$.
\end{itemize}
Note that both the eigenvalues and eigenfunctions are local to each layer and simple to obtain. The form of the eigenvalues $\lambda_{i,n}$ and normalized eigenfunctions $\widehat{\phi}_{i,n}$ for $n = 0,1,\hdots$ can be found in \cite[Appendix B]{carr_2016a}. 
The coefficients $c_{i,n}(t)$ satisfy the ordinary differential equation:
\begin{gather}
\label{eq:cin_ode}
\frac{dc_{i,n}}{dt} + D_{i}\lambda_{i,n}^{2}c_{i,n} = D_{i}\left[g_{i-1}(t)\beta_{i,3,n} + g_{i}(t)\beta_{i,4,n}\right] - \left[g_{i-1}'(t)\beta_{i,1,n} + g_{i}'(t)\beta_{i,2,n}\right],
\end{gather}
which is formulated by taking the inner product of both sides of (\ref{eq:pde_vi}) with $\widehat{\phi}_{i,n}(x)$ and setting:
\begin{gather}
\label{eq:beta_1}
\beta_{i,1,n} = \int_{l_{i-1}}^{l_{i}}\psi_{i,1}\widehat{\phi}_{i,n}\,dx,\quad
\beta_{i,2,n} = \int_{l_{i-1}}^{l_{i}}\psi_{i,2}\widehat{\phi}_{i,n}\,dx,\\ 
\label{eq:beta_3}
\beta_{i,3,n} = \int_{l_{i-1}}^{l_{i}}\psi_{i,1}''\widehat{\phi}_{i,n}\,dx,\quad
\beta_{i,4,n} = \int_{l_{i-1}}^{l_{i}}\psi_{i,2}''\widehat{\phi}_{i,n}\,dx.
\end{gather}
The general solution of (\ref{eq:cin_ode}) is given by
\begin{multline}
\label{eq:cin_general_solution}
c_{i,n}(t) = c_{i,n}(0)e^{-tD_{i}\lambda_{i,n}^{2}}\\ + D_{i}\beta_{i,3,n}\int_{0}^{t}g_{i-1}(\tau)e^{-(t-\tau)D_{i}\lambda_{i,n}^{2}}\,d\tau + D_{i}\beta_{i,4,n}\int_{0}^{t}g_{i}(\tau)e^{-(t-\tau)D_{i}\lambda_{i,n}^{2}}\,d\tau\\
-\beta_{i,1,n}\int_{0}^{t}g_{i-1}'(\tau)e^{-(t-\tau)D_{i}\lambda_{i,n}^{2}}\,d\tau
-\beta_{i,2,n}\int_{0}^{t}g_{i}'(\tau)e^{-(t-\tau)D_{i}\lambda_{i,n}^{2}}\,d\tau.
\end{multline}
The unknown coefficient $c_{i,n}(0)$ is identified using the initial condition (\ref{eq:vi_ic}):
\begin{align*}
\sum_{n=0}^{\infty} c_{i,n}(0)\widehat{\phi}_{i,n}(x) = \tilde{f}_{i}(x),
\end{align*}
and hence we have that:
\begin{gather*}
c_{i,n}(0) = \int_{l_{i-1}}^{l_{i}}\tilde{f}_{i}\widehat{\phi}_{i,n}\,dx = \int_{l_{i-1}}^{l_{i}}f_{i}\widehat{\phi}_{i,n}\,dx- g_{i-1}(0)\beta_{i,1,n} - g_{i}(0)\beta_{i,2,n}.
\end{gather*}
Applying integration by parts to the third and fourth integrals in the expression for $c_{i,n}(t)$ (\ref{eq:cin_general_solution}) yields:
\begin{multline*}
c_{i,n}(t) = \beta_{i,5,n}e^{-tD_{i}\lambda_{i,n}^{2}} - g_{i-1}(t)\beta_{i,1,n} - g_{i}(t)\beta_{i,2,n}\\ + D_{i}(\beta_{i,3,n}+\lambda_{i,n}^{2}\beta_{i,1,n})\int_{0}^{t}g_{i-1}(\tau)e^{-(t-\tau)D_{i}\lambda_{i,n}^{2}}\,d\tau\\ + D_{i}(\beta_{i,4,n}+\lambda_{i,n}^{2}\beta_{i,2,n})\int_{0}^{t}g_{i}(\tau)e^{-(t-\tau)D_{i}\lambda_{i,n}^{2}}\,d\tau,
\end{multline*}
where we have set
\begin{gather*}
\beta_{i,5,n} = \int_{l_{i-1}}^{l_{i}}f_{i}\widehat{\phi}_{i,n}\,dx.
\end{gather*}
We express the above expression for $c_{i,n}(t)$ in terms of the inverse Laplace transform for reasons that will become clear later. Using the convolution property of the Laplace transform allows an equivalent expression to be obtained:
\begin{multline}
\label{eq:cin_final}
c_{i,n}(t) = \beta_{i,5,n}e^{-tD_{i}\lambda_{i,n}^{2}} - g_{i-1}(t)\beta_{i,1,n} - g_{i}(t)\beta_{i,2,n}\\ + D_{i}(\beta_{i,3,n}+\lambda_{i,n}^{2}\beta_{i,1,n})\mathcal{L}^{-1}\left\{\frac{\overline{g}_{i-1}(s)}{s+D_{i}\lambda_{i,n}^{2}}\right\}\\ + D_{i}(\beta_{i,4,n}+\lambda_{i,n}^{2}\beta_{i,2,n})\mathcal{L}^{-1}\left\{\frac{\overline{g}_{i}(s)}{s+D_{i}\lambda_{i,n}^{2}}\right\},
\end{multline}
where, for example, $\overline{g}_{i}(s) = \mathcal{L}\{g_{i}(t)\}$. It follows therefore that our solution approach assumes the Laplace transformations of the boundary functions, namely $\overline{g}_{0}(s)$ and $\overline{g}_{m}(s)$, are known, that is, they exist and can be evaluated analytically or numerically. With $c_{i,n}(t)$ identified, we have solved (\ref{eq:v_first_layer})--(\ref{eq:v_last_layer}) for $v_{i}(x,t)$ ($i = 1,\hdots,m$). The solution of the single layer problems (\ref{eq:u_first_layer})--(\ref{eq:u_last_layer}) are hence given by: 
\begin{gather}
\label{eq:solution}
u_{i}(x,t) = g_{i-1}(t)\psi_{i,1}(x) + g_{i}(t)\psi_{i,2}(x) +\sum_{n=0}^{\infty}c_{i,n}(t)\widehat{\phi}_{i,n}(x),
\end{gather}
for all $i = 1,\hdots,m$.

\subsection{Evaluation of the solution expressions}
To evaluate the solution expressions, described by equations (\ref{eq:cin_final}) and (\ref{eq:solution}), the summations are truncated after a finite number of terms/eigenvalues $N$:
\begin{gather}
\label{eq:solution_Nterms}
{u_{i}(x,t) \simeq g_{i-1}(t)\psi_{i,1}(x)} + g_{i}(t)\psi_{i,2}(x) +\sum_{n=0}^{N-1}c_{i,n}(t)\widehat{\phi}_{i,n}(x),
\end{gather}
for all $i = 1,\hdots,m$. The inverse Laplace transformations appearing in the coefficients (\ref{eq:cin_final}) are evaluated using the strategy proposed by \citet{carr_2016a}, which involves applying the quadrature formula described by \citet{trefethen_2006} to the integral representation of the inverse Laplace transform. Given $F(s) = \mathcal{L}\{f(t)\}$, the Laplace transform can be inverted numerically as follows:
\begin{gather}
\label{eq:inverse_laplace_approx}
f(t) \approx -2\Re\left\{\sum_{k=1}^{N_{p}/2} c_{2k-1}\frac{F(z_{2k-1}/t)}{t}\right\},
\end{gather}
where $\Re\{\cdot\}$ denotes the real part, and $c_{2k-1}, z_{2k-1}\in\mathbb{C}$ {are the residues and poles of the best $(N_{p},N_{p})$ rational approximation to $e^{z}$ on the negative real line} (see, e.g., \cite{trefethen_2006,carr_2016a})\footnote{{Note that because we are dealing with diffusion problems, our interest is in functions that have the functional form $f(t) = \nu_{0} + \sum_{n=1}^{\infty} \nu_{n}e^{-\mu_{n}t}$, where $\nu_{n},\mu_{n}\in\mathbb{R}$. In this case, $F(s) = \nu_{0}/s + \sum_{n=1}^{\infty} \nu_{n}/(s-\mu_{n})$, which yields $F(z_{2k-1}/t)/t = \nu_{0}/z_{2k-1} + \sum_{n=1}^{\infty} \nu_{n}/(z_{2k-1}+\mu_{n}t)$ and therefore the limit of (\ref{eq:inverse_laplace_approx}) is well behaved as $t\rightarrow 0$ and $t\rightarrow\infty$.}}. Using (\ref{eq:inverse_laplace_approx}) gives, for example, the following approximation \cite{carr_2016a}:

\begin{gather}
\label{eq:inverse_laplace_approx_gi}
\mathcal{L}^{-1}\left\{\frac{\overline{g}_{i}(s)}{s+D_{i}\lambda_{i,n}^{2}}\right\} \approx -2\Re\left\{\sum_{k=1}^{N_{p}/2} \frac{c_{2k-1}\overline{g}_{i}(z_{2k-1}/t)}{z_{2k-1}+D_{i}\lambda_{i,n}^{2}t}\right\}.
\end{gather}
Note that the quadrature formula (\ref{eq:inverse_laplace_approx_gi}) requires evaluating $\overline{g}_{i}(s)$, which is the Laplace transformation of the unknown interface function $g_{i}(t)$. To compute these evaluations the solutions are constrained to satisfy the interface condition (\ref{eq:original_constraint}) \cite{carr_2016a}. Taking Laplace transforms of (\ref{eq:original_constraint}) and rearranging yields:
\begin{align}
\label{eq:constraint_laplace}
\overline{u}_{i}(l_{i},s) - \theta_{i}\overline{u}_{i+1}(l_{i},s) + \frac{1}{H_{i}}\overline{g}_{i}(s) = 0,
\end{align}
where $\overline{u}_{i}(x,s) = \mathcal{L}\{u_{i}(x,t)\}$ is given by:
\begin{gather}
\label{eq:laplace_tranform_ui}
{\overline{u}_{i}(x,s) \simeq \overline{g}_{i-1}(s)\psi_{i,1}(x)} + \overline{g}_{i}(s)\psi_{i,2}(x) +\sum_{n=0}^{N-1}\overline{c}_{i,n}(s)\widehat{\phi}_{i,n}(x).
\end{gather}
The Laplace transformation $\overline{c}_{i,n}(s) = \mathcal{L}\left\{c_{i,n}(t)\right\}$ is linear in the functions $g_{i}(s)$ ($i = 1,\hdots,m$)
\begin{multline}
\label{eq:cins}
\overline{c}_{i,n}(s) = \frac{\beta_{i,5,n}}{s + D_{i}\lambda_{i,n}^{2}} + \left(\frac{D_{i}(\beta_{i,3,n}+\lambda_{i,n}^{2}\beta_{i,1,n})}{s+D_{i}\lambda_{i,n}^{2}}-\beta_{i,1,n}\right)\overline{g}_{i-1}(s)\\ + \left(\frac{D_{i}(\beta_{i,4,n}+\lambda_{i,n}^{2}\beta_{i,2,n})}{s+D_{i}\lambda_{i,n}^{2}}-\beta_{i,2,n}\right)\overline{g}_{i}(s).
\end{multline}
With the forms of $\overline{g}_{0}(s)$ and $\overline{g}_{m}(s)$ given, substituting (\ref{eq:laplace_tranform_ui}) into (\ref{eq:constraint_laplace}) for all $i = 1,\hdots,m$ produces a tridiagonal matrix system in the form:
\begin{align}
\label{eq:linear_system}
\mathbf{A}(s)\mathbf{x} = \mathbf{b}(s),
\end{align}
where $\mathbf{A}(s)\in\mathbb{C}^{(m-1)\times(m-1)}$, $\mathbf{b}(s)\in\mathbb{C}^{m-1}$ and $\mathbf{x} = \left(\overline{g}_{1}(s),\hdots,\overline{g}_{m-1}(s)\right)^{T}$. The individual entries of $\mathbf{A}(s)$ and $\mathbf{b}(s)$ are given in Appendix \ref{app:linear_system}. Solving the linear system (\ref{eq:linear_system}) evaluated at $s = z_{2k-1}/t$ allows the required evaluations $\overline{g}_{i}(z_{2k-1}/t)$ ($i = 1,\hdots,m-1$) to be computed\footnote{{Note that in the limit as $t\rightarrow 0$ ($s\rightarrow\infty$) and $t\rightarrow\infty$ ($s\rightarrow 0$) the entries of $\mathbf{A}(s)$ and $\mathbf{b}(s)$, which arise from the coefficients of $\overline{g}_{i-1}(s)$ and $\overline{g}_{i}(s)$ and the first term in the expression for $\overline{c}_{i,n}(s)$ (\ref{eq:cins}) (see Appendix \ref{app:linear_system}) are well behaved.}}. For each time $t>0$ at which the semi-analytical solution is sought, the $m-1$ dimensional matrix system (\ref{eq:linear_system}) must be solved $N_{p}/2$ times.

In contrast to the case of time-independent BCs \cite{carr_2016a}, for $i=0$ and $i=m$, the coefficients (\ref{eq:cin_final}) feature terms involving inverse Laplace transformations of expressions involving $\overline{g}_{0}(s) = \mathcal{L}\{g_{0}(t)\}$ and $\overline{g}_{m}(s) = \mathcal{L}\{g_{m}(t)\}$. Since the external boundary functions $g_{0}(t)$ and $g_{m}(t)$ are known \textit{a priori}, these terms can be computed directly using the convolution property, for example:
\begin{gather}
\label{eq:laplace_inverse_g0}
\mathcal{L}^{-1}\left\{\frac{\overline{g}_{0}(s)}{s+D_{1}\lambda_{1,n}^{2}}\right\} = \int_{0}^{t} g_{0}(\tau)e^{-(t-\tau)D_{1}\lambda_{1,n}^{2}}\,d\tau.
\end{gather}
For a general and flexible code and to reduce user input, numerically evaluating the integral (\ref{eq:laplace_inverse_g0}) is preferred. However, preliminary investigation found that the quadrature formula (\ref{eq:inverse_laplace_approx_gi}) performed better {for large $t$ than MATLAB's in-built \texttt{integral} function, which uses global adaptive Gauss-Kronrod quadrature}. Therefore, we use (\ref{eq:inverse_laplace_approx_gi}) for all $i = 0,\hdots,m$.

The only thing left to address is evaluating the unknown interface functions $g_{i}(t)$ ($i=1,\hdots,m-1$) at a given time $t$. Note these evaluations appear in both the expression for $u_{i}(x,t)$ (\ref{eq:solution}) and the coefficients $c_{i,n}(t)$ (\ref{eq:cin_final}). Since $g_{i}(t) = \mathcal{L}^{-1}\left\{\overline{g}_{i}(s)\right\}$, we can use the approximation (\ref{eq:inverse_laplace_approx}), giving:
\begin{alignat}{2}
\label{eq:laplace_transformations3}
g_{i}(t) \approx -2\Re\left\{\sum_{k=1}^{N_{p}/2} \frac{c_{2k-1}\overline{g}_{i}(z_{2k-1}/t)}{t}\right\},
\end{alignat}
for $i = 1,\hdots,m-1$. Note here that the evaluations of $\overline{g}_{i}(s)$ are the same as those appearing in (\ref{eq:inverse_laplace_approx_gi}). Finally, we remark that it is due to approximations such as (\ref{eq:inverse_laplace_approx_gi}) that both the method developed in this paper and the method of \citet{carr_2016a} are classified as \textit{semi-analytical}.

\subsection{MATLAB code}
The semi-analytical solution developed originally by \citet{carr_2016a} is available at\\ \href{https://github.com/elliotcarr/MultDiff}{https://github.com/elliotcarr/MultDiff}. This repository has now been updated to include the semi-analytical solution developed in this paper, which is applicable to time-dependent external BCs and more general interface conditions.

\section{Results and discussion}
\label{sec:results}
\subsection{Effect of interface conditions}
\label{sec:interface_conditions}
\begin{figure*}[t]
\centering
\subfloat[Case A]{\includegraphics[width = 0.48\textwidth]{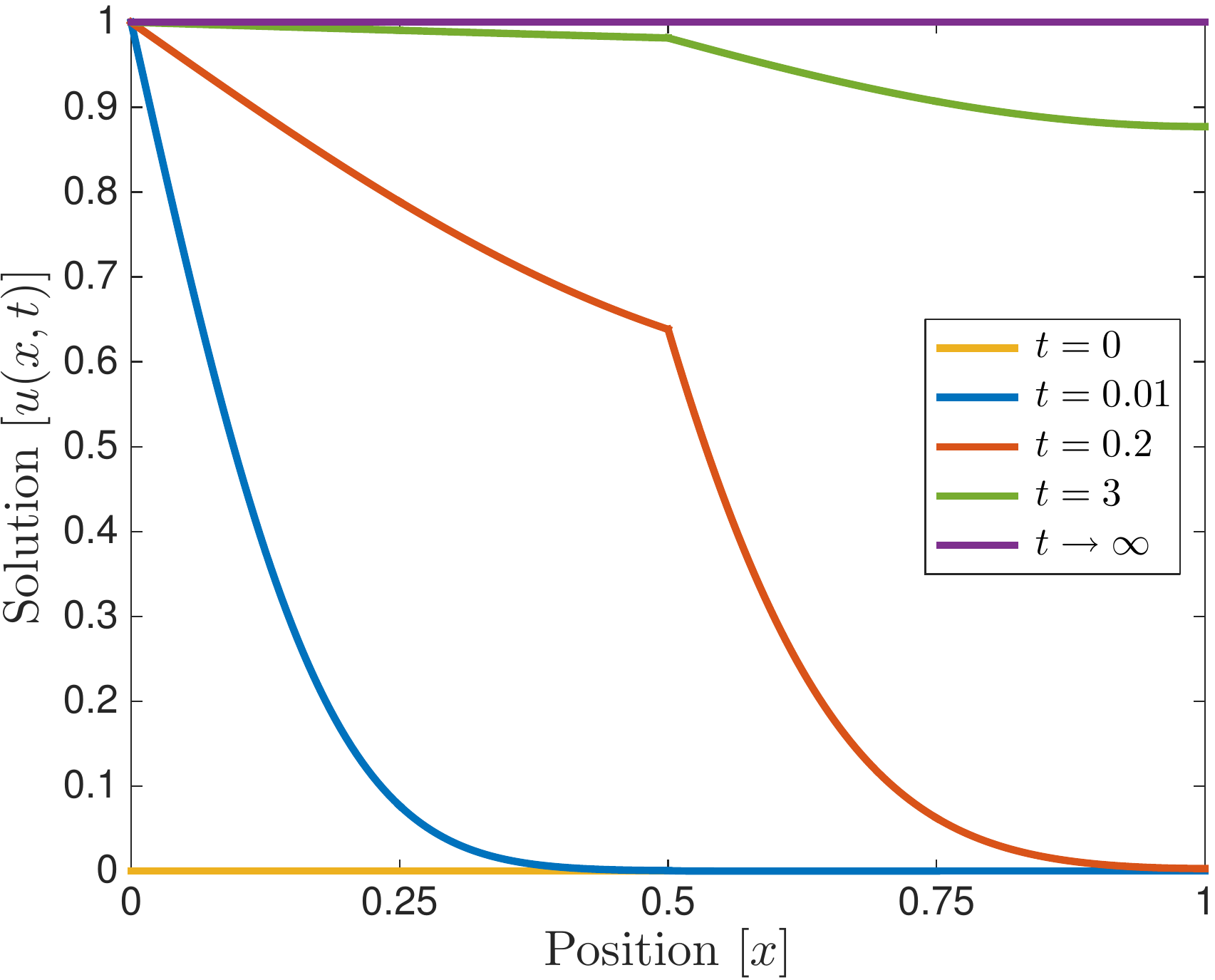}\label{fig:CaseA_interface}}\hspace{0.3cm}\subfloat[Case B]{\includegraphics[width = 0.48\textwidth]{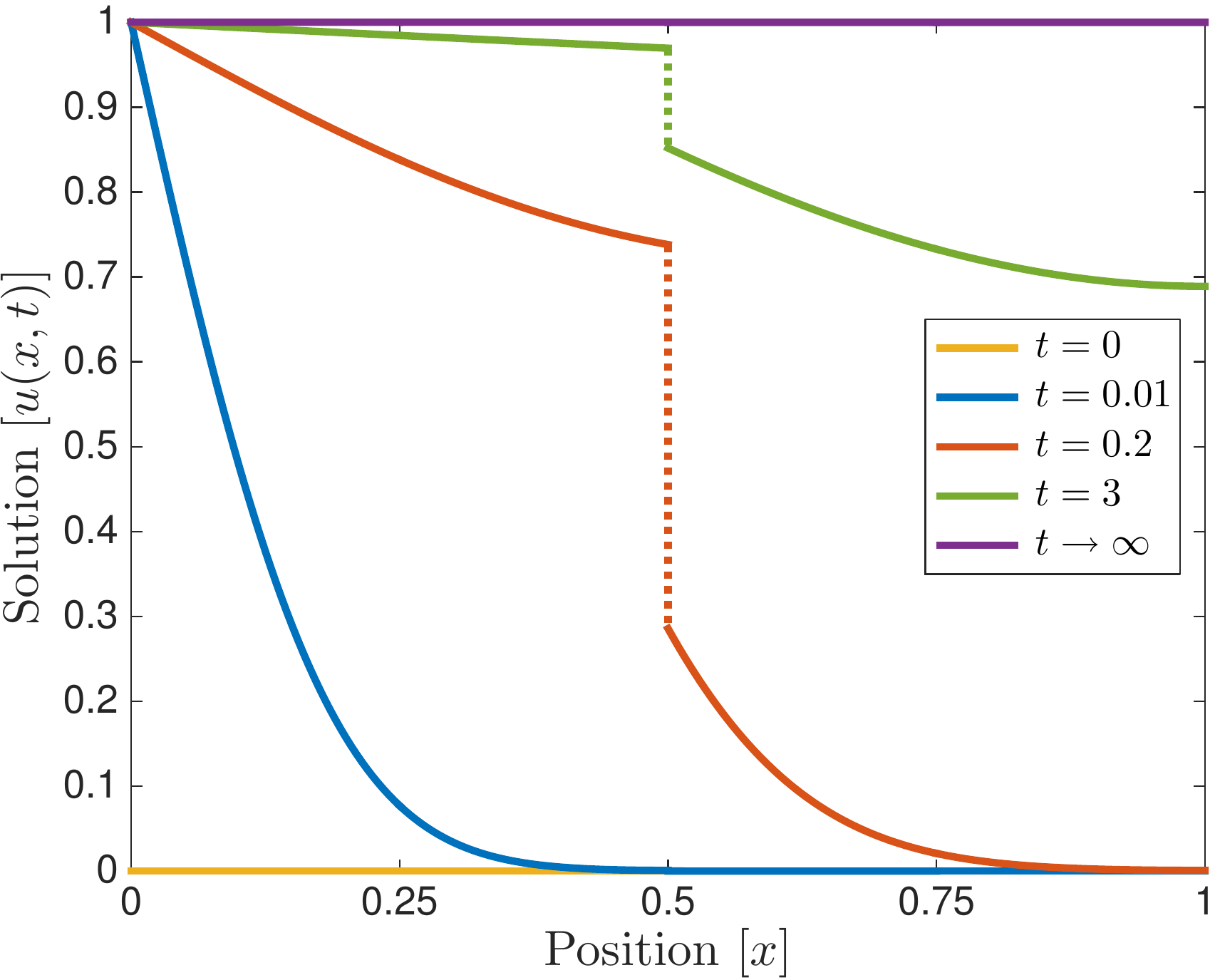}\label{fig:CaseB_interface}}\\
\subfloat[Case C]{\includegraphics[width = 0.48\textwidth]{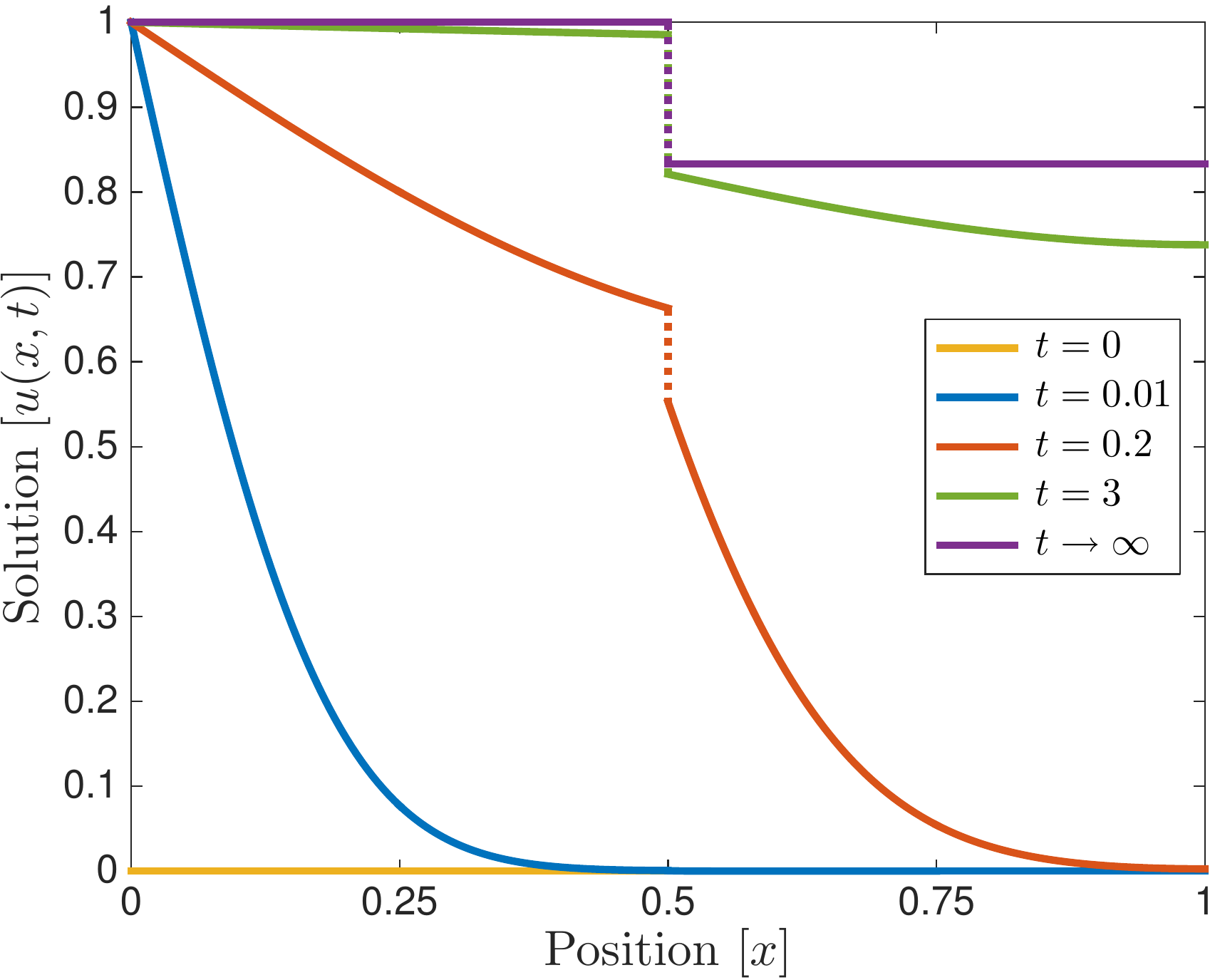}\label{fig:CaseC_interface}}\hspace{0.3cm}\subfloat[Case D]{\includegraphics[width = 0.48\textwidth]{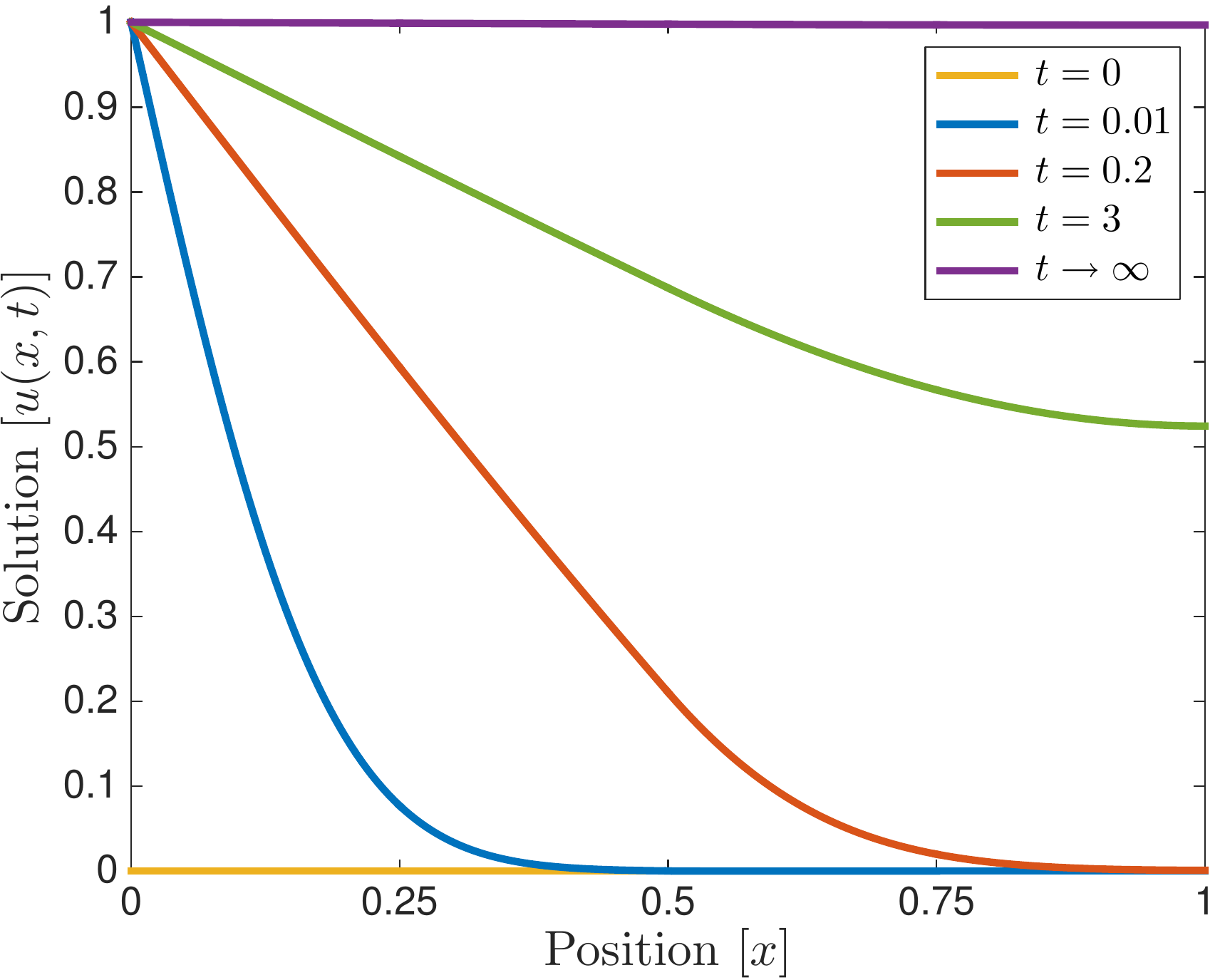}\label{fig:CaseD_interface}}
\caption{Semi-analytical solution for Cases A-D.}
\end{figure*}

Consider the multilayer diffusion problem described by equations (\ref{eq:original_equation})--(\ref{eq:original_bc2}) and (\ref{eq:ic_general}) with $m = 2$ layers, domain $[l_{0},l_{1},l_{2}] = [0,0.5,1]$, diffusivities $D_{1} = 1$ and $D_{2} = 0.1$, and external boundary condition data $a_{L} = 1$, $b_{L} = 0$, $g_{0}(t) = 1$, $a_{R} = 0$, $b_{R} = 1$ and $g_{m}(t) = 0$. Figures \ref{fig:CaseA_interface}--\ref{fig:CaseD_interface} depict the semi-analytical solution over time for four different choices of interface conditions:
\begin{itemize}
\item Case A: Perfect contact with $\gamma_{1} = D_{1}$ and $\gamma_{2} = D_{2}$.
\item Case B: Jump conditions with $\gamma_{1} = D_{1}$, $\gamma_{2} = D_{2}$ and $H_{1} = 0.5$.
\item Case C: Partition conditions with $\gamma_{1} = D_{1}$, $\gamma_{2} = D_{2}$ and $\theta_{1} = 1.2$.
\item Case D: Perfect contact with $\gamma_{1} = \gamma_{2} = 2.0$.
\end{itemize}
For Case A, the solution gradient is discontinuous at the interface (see Figure \ref{fig:CaseA_interface}) as $\gamma_{1}\neq \gamma_{2}$ (see interface condition (\ref{eq:ic_perfect_2})) while for Case B, the solution is discontinuous at the interface (see Figure \ref{fig:CaseB_interface}) since the transfer coefficient $H_{1}$ is finite. Using equation (\ref{eq:ic_jump_1}), the step change in the solution at the interface can be expressed as:
\begin{align}
\label{eq:jump_type2}
u_{2}(l_{1},t) - u_{1}(l_{1},t) = \frac{\gamma_{1}}{H_{1}}\frac{\partial u_{1}}{\partial x}(l_{1},t).
\end{align}
For given values of $H_{1}$ and $\gamma_{1}$, the difference in solution values at the interface is proportional to the gradient appearing in (\ref{eq:jump_type2}), which explains why the jump discontinuity is absent from the solution at $t = 0.01$ (at least visibly), large at $t = 0.2$ and small at $t = 5.0$ (see Figure \ref{fig:CaseB_interface}). Case C also exhibits a jump discontinuity at the interface with:
\begin{align}
\label{eq:jump_type3}
u_{2}(l_{1},t) - u_{1}(l_{1},t) = (1-\theta_{1})u_{2}(l_{1},t),
\end{align}
which means that, for a given value of the partition coefficient $\theta_{1}$, the size of the jump discontinuity is directly proportional to the value of $u_{2}(l_{1},t)$. This is confirmed in Figure \ref{fig:CaseC_interface} with the step change in the solution across the interface (\ref{eq:jump_type3}) growing over time. In contrast to the jump conditions of Case B (Figure \ref{fig:CaseB_interface}), the steady state solution is dependent upon the partition coefficient. For Case $D$, both the solution and gradient are continuous at the interface even though the diffusivities $D_{1}$ and $D_{2}$ are not equal (Figure \ref{fig:CaseD_interface}). This is explained by the interface condition  (\ref{eq:ic_perfect_2}) describing continuity of the flux, which reduces to continuity of the gradient if $\gamma_{1} = \gamma_{2}$. Indeed, because of this cancellation, the solutions are invariant under the condition $\gamma_{1} = \gamma_{2} = \gamma$ regardless of the value of $\gamma$.

\subsection{Convergence of semi-analytical solution}
\label{sec:convergence}
In this section, we investigate the rate of convergence exhibited by our semi-analytical method by comparing the error against the number of eigenvalues $N$ used in the solution expansions (\ref{eq:solution_Nterms}). Let $u_{i}(x_{i,j},t)$ denote the exact solution of the multilayer diffusion problem in the $i$th layer evaluated at grid point $x_{i,j} = l_{i-1}+j\Delta x_{i}$ where the grid spacing $\Delta x_{i} := (l_{i}-l_{i-1})/N_{x}$ and $N_{x}$ is the constant number of divisions in each layer. Furthermore, let $\widehat{u}_{i}^{(N)}(x_{i,j},t)$ denote an approximate analytical solution computed using $N$ terms/eigenvalues in each layer. The relative error of this approximate analytical solution is computed as
\begin{align}
\label{eq:relative_error}
\varepsilon_{N}(t) = \frac{\max\limits_{i,j} \left|u_{i}(x_{i,j},t) - \widehat{u}_{i}^{(N)}(x_{i,j},t)\right|}{\max\limits_{i,j} \left|u_{i}(x_{i,j},t)\right|},
\end{align}
where the maximum is taken over $i = 1,\hdots,m$ and $j = 1,\hdots,N_{x}+1$. Using the above error definition, we compare the accuracy of the semi-analytical solution (\ref{eq:solution_Nterms}) to the classical analytical solution derived using separation of variables (see, e.g., \cite{carr_2016a,hickson_2009a,trefry_1999}). The classical analytical solution subject to the general interface conditions (\ref{eq:ic_general}) is described briefly in Appendix \ref{app:analytical_solution}.

\begin{figure}[h]
\centering
\subfloat[$D_{1} = 1$ and $D_{2} = 0.1$]{\label{fig:error_small}\includegraphics[width = 0.49\textwidth]{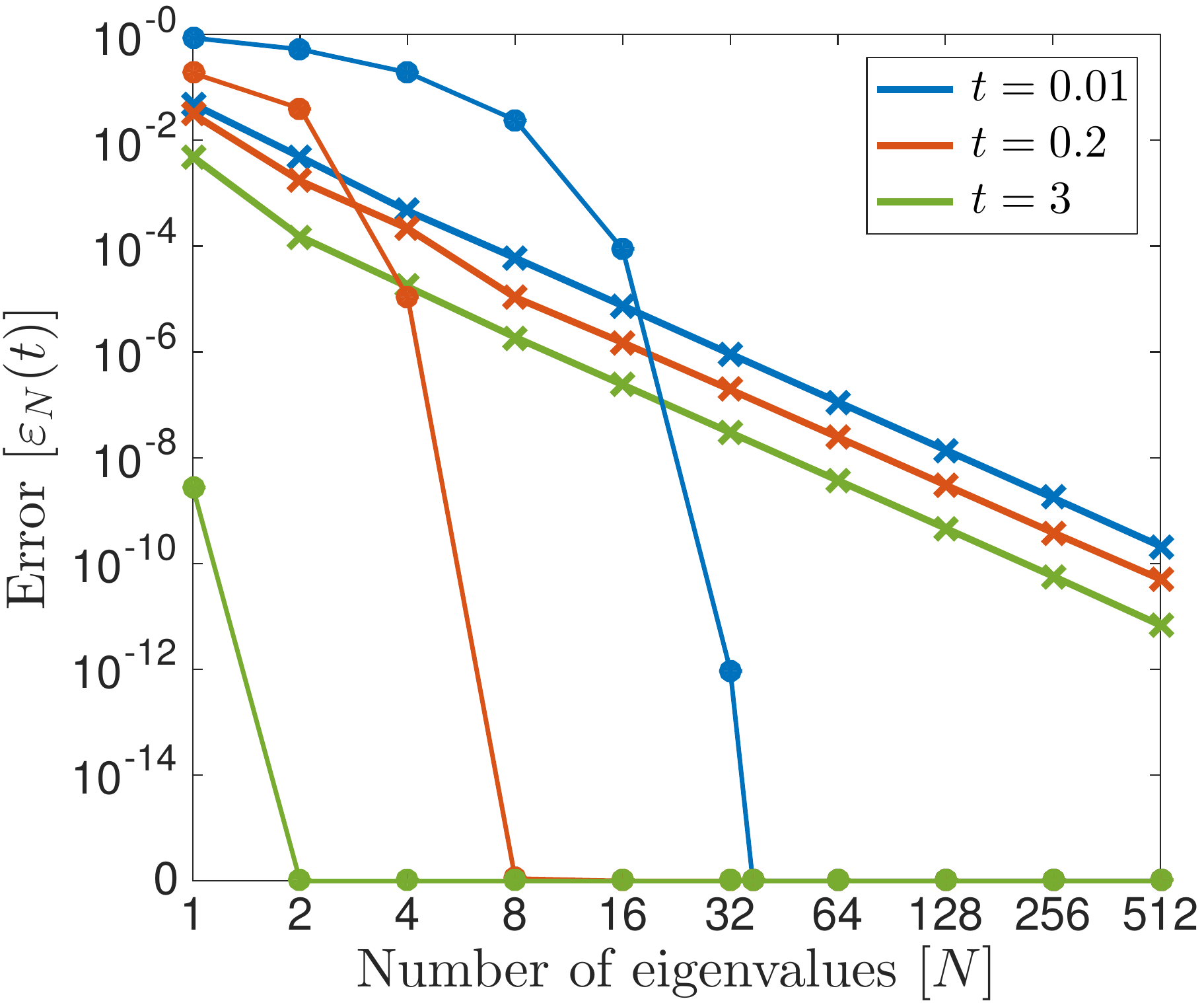}}
\subfloat[$D_{1} = 100$ and $D_{2} = 0.1$]{\label{fig:error_large}\includegraphics[width = 0.49\textwidth]{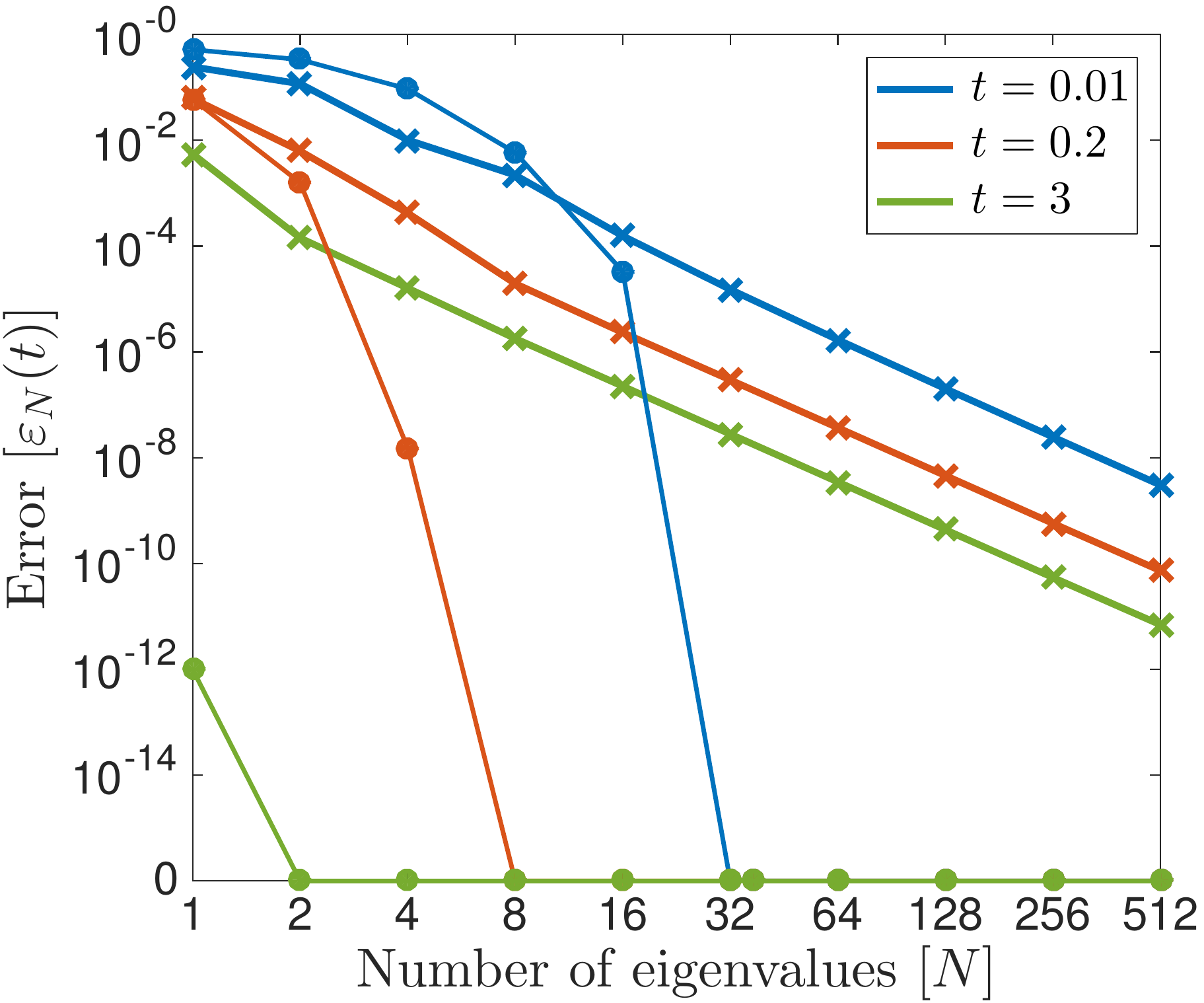}}
\caption{Plot of the error [$\varepsilon_{N}(t)$] for the semi-analytical solution (crosses) and analytical solution (dots) (computed using $N_{x}=5$) versus the number of eigenvalues [$N$] for Case C of Section \ref{sec:interface_conditions} and both a (a) small and (b) large contrast in diffusivity. Note the slight abuse of scale along the vertical axis depicting the value of zero on a log scale.}
\label{fig:error}
\end{figure}

Figure \ref{fig:error} gives the relative error for both the semi-analytical and analytical solutions versus the number of eigenvalues for Case $C$ of Section \ref{sec:interface_conditions} for both a small and large contrast in diffusivity. For increasing $N$, due to the presence of the exponential $e^{-\lambda_{n}^{2}t}$ in the analytical solution (see Appendix \ref{app:analytical_solution}), the terms in the solution expansion tend to zero extremely rapidly provided $t$ isn't too small: for $t = 0.01,0.2,3$ the exact solution is effectively obtained (i.e., the error falls below the machine epsilon of $2^{-52}$ in MATLAB) after $N = 37,8,2$ eigenvalues, respectively. 

This behaviour is not observed for the semi-analytical solution and this is the price paid for reformulating the problem (see Section \ref{sec:reformulation}) to avoid solving a complex transcendental equation arising from a matrix determinant to determine the eigenvalues. While the first term of $c_{i,n}(t)$ (\ref{eq:cin_final}) approaches zero extremely rapidly for increasing $N$, the second two terms, which arise due to the reformulation, do not. Little difference is observed when comparing Figures \ref{fig:error_small} and \ref{fig:error_large}, which demonstrates that the semi-analytical method performs well for both small and large contrasts in the diffusivity.

Both plots in Figure \ref{fig:error} suggest the following linear relationship exists between the logarithm of the error $\varepsilon_{N}(t)$ for the semi-analytical method and the number of eigenvalues $N$:
\begin{gather*}
\log\varepsilon_{N}(t) \approx \log\varepsilon_{1}(t) - p\log N,
\end{gather*}
or equivalently:
\begin{gather*}
\varepsilon_{N}(t) \approx \varepsilon_{1}(t)N^{-p},
\end{gather*}
where $p$ denotes the negative slope of the linear error curves given in Figures \ref{fig:error_small} and \ref{fig:error_large}, which is essentially independent across the three values of time $t$. Successive computation of the slopes for increasing $N$ indicates that for each time $t$, the slope $p$ approximately approaches a value of $3$ as $N$ increases. Together these observations suggest the following {proposition. 
\medskip
\begin{proposition} The convergence rate for the semi-analytical method is
\begin{gather}
\label{eq:converence_rate}
\varepsilon_{N}(t) = O(N^{-3})\quad \text{for large $N$}.
\end{gather}
\end{proposition}}
\begin{proof}
\textit{{This result can be derived} by studying the coefficient $c_{i,n}(t)$ appearing in the solution expression (\ref{eq:solution_Nterms}) and noting that each integral in equation (\ref{eq:cin_general_solution}) is $O(\lambda_{i,n}^{-2})$ for large $\lambda_{i,n}$. For example, applying integration by parts to the second integral in equation (\ref{eq:cin_general_solution}) yields:
\begin{align}
\label{eq:integral_order}
\int_0^t g_i(\tau) e^{-(t-\tau)D_i\lambda_{i,n}^2} d\tau = \frac{g_i(t) - g_i(0)e^{-tD_i\lambda_{i,n}^2}}{D_{i}\lambda_{i,n}^2} + O\left(\lambda_{i,n}^{-4}\right) = O\left(\lambda_{i,n}^{-2}\right),
\end{align}
for large $\lambda_{i,n}$. Additionally, one can show that each of the constants defined in equation (\ref{eq:beta_1}) and (\ref{eq:beta_3}) is either $O(\lambda_{i,n}^{-1})$ or equal to zero. For example, consider the values of $\beta_{i,1,n}$  (\ref{eq:beta_1}) and $\beta_{i,3,n}$ (\ref{eq:beta_3}) in the first layer ($i=1$) subject to a Dirichlet BC at $x = l_{0}$. Since: 
\begin{align*}
\psi_{1,1}(x) = \frac{1}{a_{L}},\quad
\widehat{\phi}_{1,n}(x) = \sqrt{\frac{2}{l_{1}-l_{0}}}\sin(\lambda_{1,n}(x-l_{0})),
\end{align*}
(see Appendix \ref{app:psi_functions} and Appendix B of \citet{carr_2016a}, respectively), clearly $\beta_{1,3,n} = 0$ as $\psi_{1,1}''(x) = 0$ while:
\begin{align*}
\beta_{1,1,n} = -\sqrt{\frac{2}{l_{1}-l_{0}}}\frac{1}{a_{L}\lambda_{1,n}} = O(\lambda_{1,n}^{-1}),
\end{align*}
for large $\lambda_{i,n}$. Verifying the remaining cases follows similarly. With the above results $c_{i,n}(t) = O(\lambda_{i,n}^{-3})$ for large $\lambda_{i,n}$. Hence, subtracting (\ref{eq:solution_Nterms}) from (\ref{eq:solution}) gives the following expression for the error:
\begin{align}
\label{eq:error_proof}
\varepsilon_{N}(t) = \sum_{n=N}^{\infty}c_{i,n}(t)\widehat{\phi}_{i,n}(x) = O(\lambda_{i,N}^{-3}),
\end{align}
for large $\lambda_{i,N}$. By considering the possible cases for the eigenvalues (see Appendix B of \citet{carr_2016a}), we observe that $\lambda_{i,n} = O(n)$ for each case involving explicit expressions for the eigenvalues. This is also true for the remaining cases (i.e., the first and last layers under Robin BCs) since the eigenvalues, $\lambda_{1,n}$ and $\lambda_{m,n}$, get closer and closer to $n\pi/(l_{1}-l_{0})$ and $n\pi/(l_{m}-l_{m-1})$, respectively, for large $n$ (large $\lambda_{i,n}$) \cite{strauss_1992}. It follows that $\lambda_{i,N}^{-3} = O(N^{-3})$ and therefore from (\ref{eq:error_proof}) we get the desired result.}\hfill$\blacksquare$
\end{proof}

We remark that this slowed convergence is typical of problems with time-dependent BCs. For example, consider the following single-layer problem with a time-dependent Dirichlet BC at $x=0$:
\begin{gather*}
\frac{\partial u}{\partial t} = D\frac{\partial^{2}u}{\partial x^{2}},\quad x \in (0,1),\quad t > 0,\\
u(x,0) = f(x),\quad u(0,t) = g_{0}(t),\quad \frac{\partial u}{\partial x}(1,t) = 0.
\end{gather*}
The solution, truncated after $N$ terms/eigenvalues, is given by
\begin{gather*}
u(x,t) = g_{0}(t) + 2\sum_{n=0}^{N-1}c_{n}(t)\sin(\lambda_{n}x),
\end{gather*}
where the eigenvalues $\lambda_{n} = (2n+1)\pi/2$ and the coefficients are defined as
\begin{gather*}
c_{n}(t) = e^{-D\lambda_{n}^{2}t}\int_{0}^{1} \big[f(x)-g_{0}(0)\big]\sin(\lambda_{n}x)\,dx - \frac{1}{\lambda_{n}}\int_{0}^{t}g_{0}'(\tau)e^{-(t-\tau)D\lambda_{n}^{2}}\,d\tau.
\end{gather*}
For this problem the convergence rate is also $O(N^{-3})$ since $c_{n}(t) = O(\lambda_{n}^{-3}) = O(n^{-3})$ due to the presence of the additional integral term involving the boundary function $g_{0}(t)$, which is $O(\lambda_{n}^{-3})$ following (\ref{eq:integral_order}).

\subsection{Comparison to Sheils' UTM Heat Code}
\label{sec:compare_Sheils}
We now compare our MATLAB implementation of the semi-analytical solution derived in Section \ref{sec:semi-analytical}, available at \href{https://github.com/elliotcarr/MultDiff}{https://github.com/elliotcarr/MultDiff}, to the unified transform method of \citet{sheils_2016}, available at \href{https://github.com/nsheils/UTM_Heat}{https://github.com/nsheils/UTM\_Heat}. The chosen test case is an $m=8$ layer version of Case A (from Section \ref{sec:interface_conditions}) with domain $[l_{0},l_{m}] = [0,1]$, interfaces $l_{i} = i/m$ for $i = 1,\hdots,m-1$, diffusivities $D_{2i-1} = 1$ and $D_{2i} = 0.1$ for $i = 1,\hdots,m/2$, and external boundary condition data $a_{L} = 1$, $b_{L} = 0$, $g_{0}(t) = 1$, $a_{R} = 0$, $b_{R} = 1$ and $g_{m}(t) = 0$. Recall that Sheils' method is faulty at the end points whenever nonhomogeneous external BCs are applied \cite{sheils_2016}. To circumvent this issue, we first decompose the solution into its steady state and transient parts: $u_{i}(x,t) = 1 + v_{i}(x,t)$ for each layer $i$, and solve for $v_{i}(x,t)$, which satisfies homogeneous external BCs. 

Table \ref{tab:compare_UTM} compares the runtimes\footnote{All tests cases were carried out in MATLAB R2014b on a MacBook Pro (mid 2014) running MAC OS X Version 10.10.5 with 16 GB of RAM and a 3.0GHz dual-core Intel Core i7 processor.} and relative errors of our semi-analytical method (for different numbers of eigenvalues $N$) to Sheils' unified transform method (with the default solver options). While the semi-analytical method is less accurate for our default number of eigenvalues ($N=50$), as reported by \citet{sheils_2016}, highly accurate solutions can be computed by simply taking more terms/eigenvalues in the solution expansions (\ref{eq:solution_Nterms}): with $N = 600$ terms, the semi-analytical solution is more accurate for all three times reported and twice as fast as the unified transform method. Moreover, for smaller values of $N$, the semi-analytical solution is very fast with an accuracy that is probably sufficient for most applications.

\begin{figure}[H]
\centering
\includegraphics[width = 0.48\textwidth]{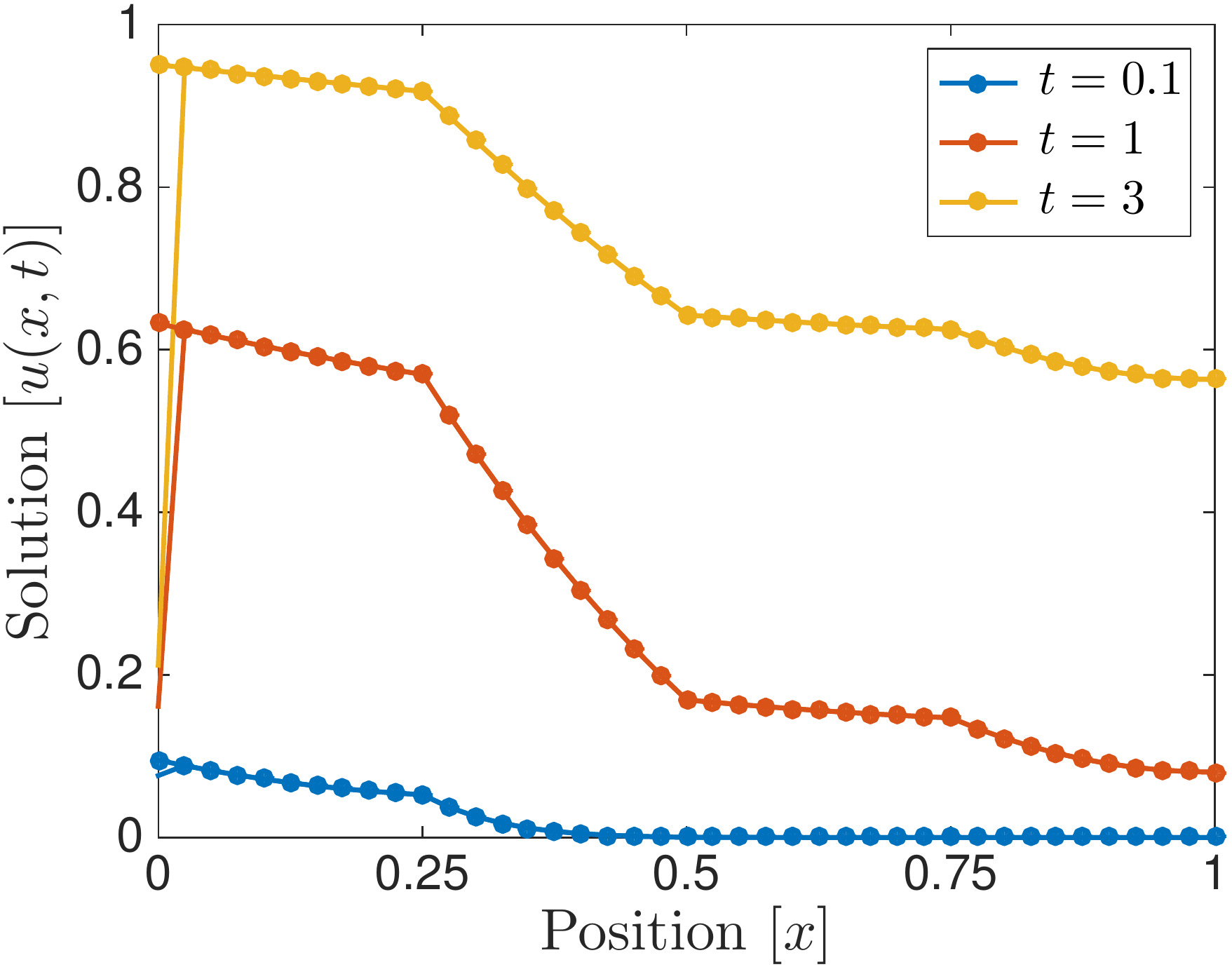}
\caption{Semi-analytical solution (dots) and unified transform method solution (continuous line) for the test case discussed in Section \ref{sec:compare_Sheils}. The unified transform method is faulty at the left boundary ($x = 0$).}
\label{fig:Faulty}
\end{figure}

For time-independent external BCs, one can always homogenise the external BCs before applying the unified transform method (as above), however, it is not clear how to do this for time-dependent external BCs. For example, with $g_{0}(t) = 1-e^{-t}$, it is not immediately obvious how to avoid the faulty behaviour of Sheils' unified transform method at the left end point (Figure \ref{fig:Faulty}). This detail, together with the fact that the analytical solution (Appendix \ref{app:analytical_solution}) doesn't perform well for a large number of layers (as reported by \citet{carr_2016a}), leads us to the conclusion that only the semi-analytical method introduced in this paper is able to correctly handle both time-dependent external BCs and a large number of layers.

\begin{table}[H]
\centering
\begin{tabular}{lrrrr}
\hline
& Runtime  & \multicolumn{3}{c}{Relative errors}\\
Method & (secs) & $t = 0.01$ & $t = 0.2$ & $t = 3$\\
\hline
Unified transform & 114.93 & 2.86e-10 & 4.65e-11 & 1.98e-11\\
Semi-analytical [$N = 10$] & 0.54 & 7.18e-05 & 2.26e-06 & 7.50e-07\\
Semi-analytical [$N = 25$] & 1.06 & 4.56e-06 & 1.43e-07 & 1.99e-08\\
Semi-analytical [$N = 50$] & 2.43 & 4.43e-07 & 1.31e-08 & 5.86e-09\\
Semi-analytical [$N = 100$] & 7.33 & 5.34e-08 & 1.70e-09 & 7.33e-10\\
Semi-analytical [$N = 200$] & 15.78 & 6.55e-09 & 1.98e-10 & 9.09e-11\\
Semi-analytical [$N = 300$] & 26.66 & 1.93e-09 & 5.98e-11 & 2.70e-11\\
Semi-analytical [$N = 400$] & 32.88 & 8.13e-10 & 2.56e-11 & 1.16e-11\\
Semi-analytical [$N = 500$] & 44.38 & 4.15e-10 & 1.29e-11 & 6.04e-12\\
Semi-analytical [$N = 600$] & 57.02 & 2.40e-10 & 7.53e-12 & 3.21e-12\\
\hline
\end{tabular}
\caption{Relative errors for our semi-analytical solution (Section \ref{sec:semi-analytical}) and the unified transform method \cite{sheils_2016} (computed using the default value of $N_{x}=15$ as per Sheils' \cite{sheils_2016} UTM Heat code) for the test case described in Section \ref{sec:compare_Sheils}.}
\label{tab:compare_UTM}
\end{table}

\subsection{Applications}
\label{sec:applications}
In this section, we apply the semi-analytical solution described in Section \ref{sec:semi-analytical} to some environmental, industrial and biological applications. The primary aim is to highlight the wide array of problems that can be solved by the code and confirm the validity of the derived semi-analytical solution by reproducing results previously reported in the literature.

\subsubsection{Contaminant transport in an aquitard}
\label{sec:contaminant_transport}
\citet{liu_1998} define the following governing equations for diffusion of a dissolved contaminant in a layered porous medium:
\begin{gather}
R_{i}\frac{\partial C_{i}}{\partial t} = D_{i}\frac{\partial^{2} C_{i}}{\partial x^{2}},
\end{gather}
for $x\in(l_{i-1},l_{i})$, $t>0$ and $i = 1,\hdots,m$, where $C_{i}(x,t)$ is the volume-based aqueous contaminant concentration [$\mu\mathrm{g}/\mathrm{L}$] in the $i$th layer, and $R_{i}$ and $D_{i}$ are constants defined as the dimensionless retardation factor and effective diffusion coefficient in the $i$th layer, respectively. At the interface between adjacent layers ($x = l_{i}$), continuity of mass flux and aqueous concentration is imposed:
\begin{gather}
C_{i}(l_{i},t) = C_{i+1}(l_{i},t),\\
\epsilon_{i}D_{i}\frac{\partial C_{i}}{\partial x}(l_{i},t) = \epsilon_{i+1}D_{i+1}\frac{\partial C_{i+1}}{\partial x}(l_{i},t),
\end{gather}
for $t > 0$ and $i = 1,\hdots,m-1$, where the constant $\epsilon_{i}$ is the porosity in the $i$th layer. 

We consider the two-layer test problem described by \citet{liu_1998}, where initially $C_{i}(x,0) = 0$ in both layers ($i = 1,2$). The concentration at the top of the first layer is assumed to be a known but arbitrary function of time, denoted by $f(t)$, while zero mass flux is assumed at the bottom boundary:
\begin{gather}
C_{1}(l_0,t) = f(t)\,,\quad\frac{\partial C_{2}}{\partial x}(l_2,t) = 0,
\end{gather}
for $t > 0$. The function $f(t)$ is assumed to take a Gaussian form:\begin{gather}
f(t) = C_{\mathrm{max}}\exp\left(-\frac{(t-\mu)^2}{\sigma^{2}}\right),
\end{gather}
where the constant $C_{\mathrm{max}}$ is the peak concentration, and $\mu$ and $\sigma$ are constants. The semi-analytical solution requires the Laplace transformation of $f(t)$, which is found to be:
\begin{gather}
\overline{f}(s) = \frac{\sigma\sqrt{\pi}C_{\mathrm{max}}}{2}\left[1+\text{erf}\left(\frac{2\mu-s\sigma^{2}}{2\sigma}\right)\right]\exp\left(\frac{s}{4}\left(s\sigma^{2}-4\mu\right)\right).
\end{gather}
Note that the quadrature formula (\ref{eq:inverse_laplace_approx_gi}) requires evaluation of $\overline{g}_{0}(s) = \overline{f}(s)$ for $s\in\mathbb{C}$. To compute the error function for complex arguments, which isn't available for the inbuilt MATLAB function \texttt{erf}, we use the \texttt{erfz} function developed by \citet{erfz_2003}. 

\begin{figure}[htb]
\centering
\includegraphics[width = 0.55\textwidth]{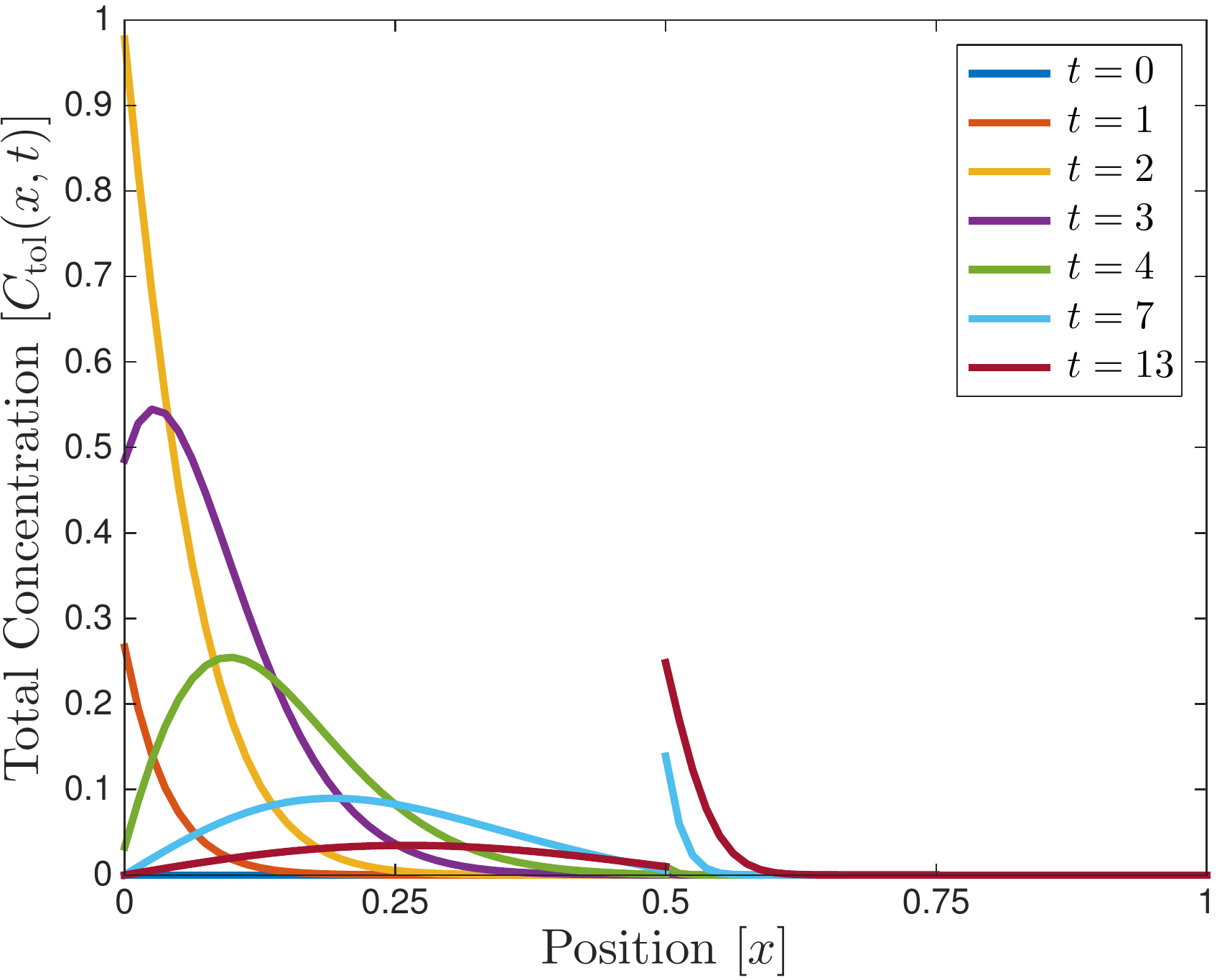}
\caption{Total concentration profiles $C_{\mathrm{tol}}(x,t)$ for the contaminant transport problem described in Section \ref{sec:contaminant_transport}. Units: $C_{tol}(x,t)$ ($\mu\mathrm{g}/\mathrm{kg}$), $x$ ($\mathrm{m}$) and $t$ ($\mathrm{yr}$).}
\label{fig:Case_Liu1}
\end{figure}

Figure \ref{fig:Case_Liu1} depicts the concentration profiles computed using the semi-analytical solution over time for $C_{\mathrm{max}} = 1.0\,\mu\mathrm{g}/\mathrm{L}$, $\mu = 2.15\,\mathrm{yr}$, $\sigma = 1\,\mathrm{yr}$, $R_{1} = 42.42$, $R_{2} = 1.67$, $D_{1} = 1.6\times 10^{-10}\,\mathrm{m}^{2}/\mathrm{s}$, $D_{2} = 2.13\times 10^{-10}\,\mathrm{m}^{2}/\mathrm{s}$ and $\epsilon_{1} = \epsilon_{2} = 0.54$. To allow comparison with the figure presented by \citet{liu_1998} the total concentration is shown in Figure \ref{fig:Case_Liu1}, which is defined in the $i$th layer as $\widetilde{C}_{i}(x,t) = C_{i}(x,t)\epsilon_{i}R_{i}/{\rho_b}_{i}$ [$\mu\mathrm{g}/\mathrm{kg}$], where ${\rho_{b}}_{i}$ is the bulk density of the soil in the $i$th layer with ${\rho_{b}}_{1} = {\rho_{b}}_{2} = 1.4\,\mathrm{kg}/\mathrm{L}$.

The plot clearly depicts the time-varying inlet concentration, reaching a peak concentration after $2.15\,\mathrm{yrs}$ ($t = \mu$). Note that the total concentration is discontinuous at the interface due to the different retardation factors in both layers. The solutions are in excellent agreement with those given by \citet{liu_1998}, giving us confidence that our new semi-analytical method has been formulated correctly for time-dependent external BCs.

\subsubsection{Heat conduction in composite materials}
\label{sec:heat_conduction}
Consider the classical heat conduction problem in a composite medium comprising $m$ layers \cite{mikhailov_1983}, where the temperature distribution in each layer is governed by the heat equation:
\begin{gather}
\label{eq:heat_equation}
\rho_{i}{c_{p}}_{i}\frac{\partial T_{i}}{\partial t} = k_{i}\frac{\partial^{2}T_{i}}{\partial x^{2}},
\end{gather}
for $x\in(l_{i-1},l_{i})$, $t > 0$ and $i = 1,\hdots,m$, and interface conditions
\begin{gather}
k_{i}\frac{\partial T_{i}}{\partial x}(l_{i},t) = H_{i}(T_{i+1}(l_{i},t) - T_{i}(l_{i},t)),\\
\label{eq:heat_ic_2}
k_{i}\frac{\partial T_{i}}{\partial x}(l_{i},t) = k_{i+1}\frac{\partial T_{i+1}}{\partial x}(l_{i},t),
\end{gather}
for $t >0$ and $i = 1,\hdots,m-1$, where $T_{i}(x,t)$ is the temperature at position $x$ and time $t$ in the $i$th layer and $H_{i}$ is the heat transfer coefficient between layers $i$ and $i+1$. The remaining constants $\rho_{i}$, ${c_{p}}_{i}$ and $k_{i}$ denote the density, specific heat capacity and thermal conductivity, respectively, in the $i$th layer. Note that in general the solution to the above problem cannot be obtained using the semi-analytical method described by \citet{carr_2016a} since it is the thermal conductivities $k_{i}$ not the thermal diffusivities $k_{i}/(\rho_{i} {c_{p}}_{i})$ that appear in the second interface condition (\ref{eq:heat_ic_2}). 

\begin{figure}[htb]
\centering
\includegraphics[width = 0.55\textwidth]{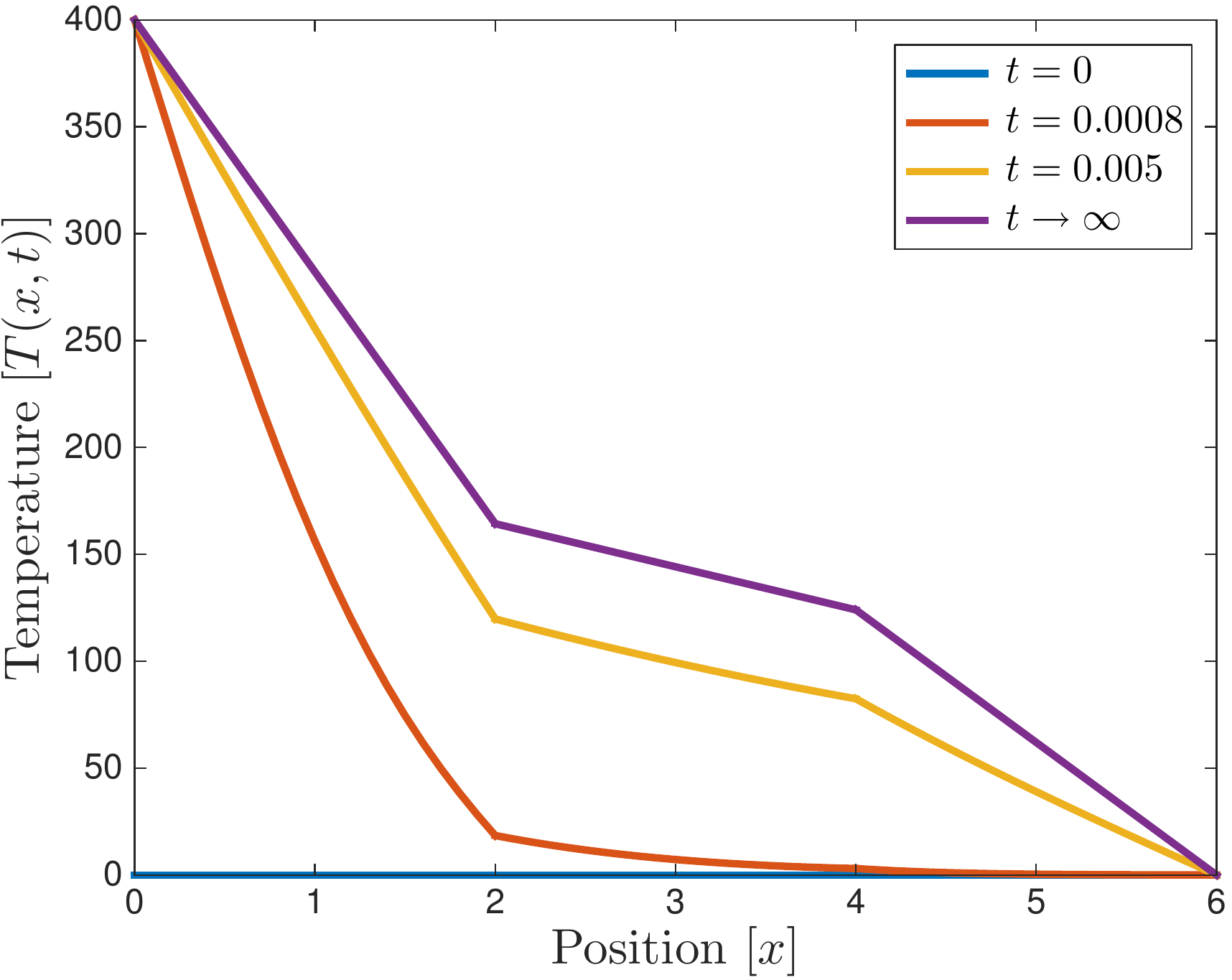}
\caption{Temperature distribution $T(x,t)$ for the heat conduction problem described in Section \ref{sec:heat_conduction}. Units: $T(x,t)$ ($^{\circ}\mathrm{C}$), $x$ ($\mathrm{cm}$) and $t$ ($\mathrm{h}$).}
\label{fig:mikhailov}
\end{figure}

As an illustrative example, we consider a classical heat conduction problem that has been solved by several authors, including \citet{mulholland_1972}, \citet{mikhailov_1983} and \citet{johnston_1991}. The test case considers $m = 3$ layers, domain $[l_{0},l_{1},l_{2},l_{3}] = [0,2,4,6]\,\mathrm{cm}$ and the following parameter values: $H_{i}\rightarrow\infty$ for all $i = 1,\hdots,3$ (i.e., perfect thermal contact at the interfaces), thermal conductivities $k_1 = 297.64$, $k_2 = 1741.18$ and $k_3 = 565.51$ [units $\mathrm{cal}/(\mathrm{cm}\,^{\circ}\mathrm{C}\,\mathrm{h}$)], densities $\rho_1 = 11.08$, $\rho_2 = 2.71$ and $\rho_3 = 7.4$ [units $\mathrm{g}/\mathrm{cm}^{3}$], and specific heat capacities $c_{p,1} = 0.031$, $c_{p,2} = 0.181$ and $c_{p,3} = 0.054$ [units $\mathrm{cal}/(\mathrm{g}\,^{\circ}\mathrm{C})$]. The initial conditions and external BCs are defined as
\begin{gather*}
T_{i}(x,0) = 0\,^{\circ}\mathrm{C},\quad i = 1,\hdots,3,\\
T_{1}(0,t) = 400\,^{\circ}\mathrm{C},\quad T_{3}(6,t) = 0\,^{\circ}\mathrm{C}.
\end{gather*}
Figure \ref{fig:mikhailov}, which shows the temperature distribution obtained using the semi-analytical solution at various times, reproduces the results previously reported in the literature \cite{mulholland_1972,mikhailov_1983,johnston_1991}.

\subsubsection{Analyte transport in composite media}
\label{sec:analyte_transport} 
Chemical concentration profiles in composite media can be sharply discontinuous at material interfaces due to partitioning phenomena \cite{todo_2013,trefry_1999}. \citet{trefry_1999} present the following mathematical model governing analyte transport in a composite medium comprising $m$ laminates/media, consisting of the diffusion equation
\begin{gather*}
\frac{\partial C_{i}}{\partial t} = \frac{D_{i}}{1+\sigma_{i}}\frac{\partial^{2} C_{i}}{\partial x^{2}},
\end{gather*}
for $x\in(l_{i-1},l_{i})$, $t > 0$ and $i = 1,\hdots,m$, subject to the following interface conditions:
\begin{gather*}
C_{i}(l_{i},t) = \alpha_{i}C_{i+1}(l_{i},t),\\ 
D_{i}\frac{\partial C_{i}}{\partial x}(l_{i},t) = D_{i+1}\frac{\partial C_{i+1}}{\partial x}(l_{i},t),
\end{gather*}
for $t > 0$ and $i =1,\hdots,m-1$. In the above equations, $C_{i}(x,t)$ is the mobile phase concentration at position $x$ and time $t$ in the $i$th layer, $\alpha_{i}$ is the mobile phase partition coefficient at the interface between layers $i$ and $i+1$ ($i = 1,\hdots,m-1$), and $D_{i}$ and $\sigma_{i}$ are the diffusion coefficient for the analyte and the linear sorption coefficient in medium $i$, respectively. Note that the solution of this problem can only be computed using the semi-analytical solution given by \citet{carr_2016a} if $\alpha_{i} = 1$ for all $i = 1,\hdots,m-1$ and $\sigma_{i} = 0$ for all $i = 1,\hdots,m$.

\begin{figure}[htb]
\centering
\includegraphics[width = 0.47\textwidth]{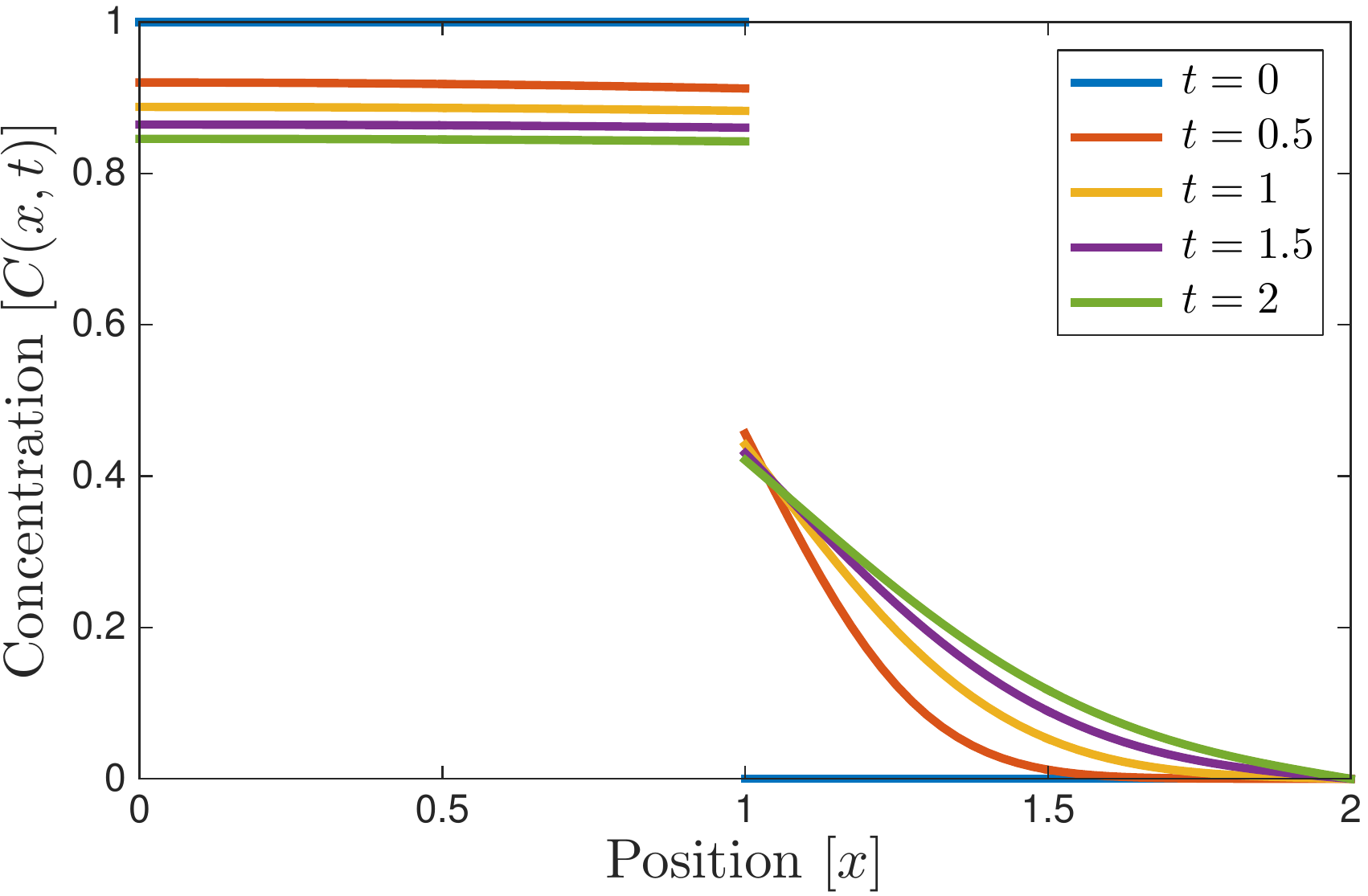}\hspace{0.3cm}\includegraphics[width = 0.47\textwidth]{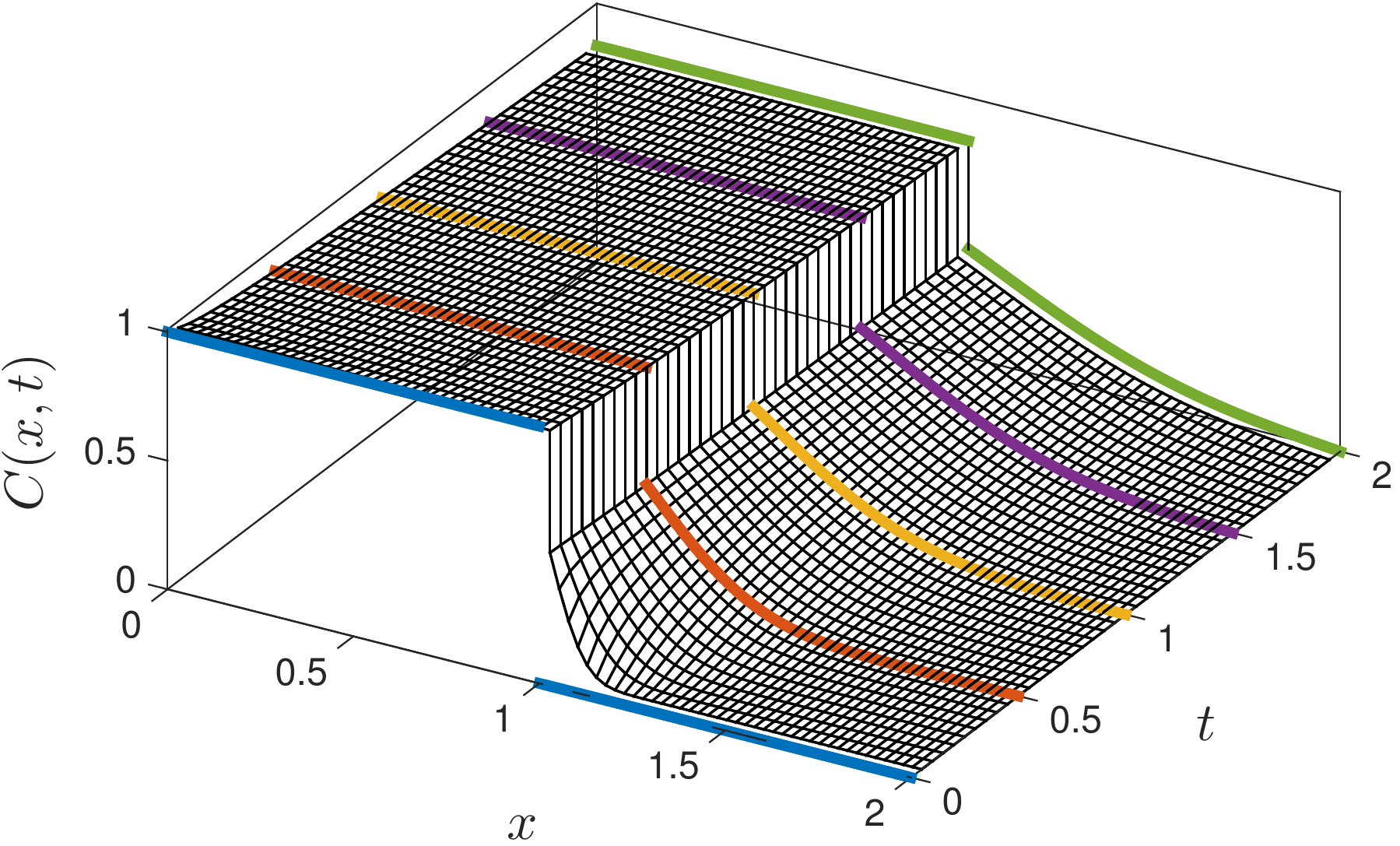}
\caption{Concentration profiles $C(x,t)$ for the analyte transport problem described in Section \ref{sec:analyte_transport}.}
\label{fig:trefry}
\end{figure}

Consider the test case described by \citet{trefry_1999} with $m=2$ layers, domain $[l_{0},l_{1},l_{2}] = [0,1,2]$, diffusion coefficients $D_{1} = 5$ and $D_{2} = 0.05$, partition coefficient $\alpha_{1} = 2$ and zero sorption coefficients $\sigma_{1}=\sigma_{2} =0$. Initially, medium $1$ is fully concentrated with analyte, medium 2 has zero initial concentration. A zero flux condition applies at the left boundary and a zero concentration condition applies at the right boundary.  

Figure \ref{fig:trefry} plots the concentration profile for $0 \leq t \leq 2$ obtained using the semi-analytical solution.  The included three-dimensional plot is in excellent agreement to the one featuring in the paper by \citet{trefry_1999}. 

\subsubsection{Brain tumour growth}
\label{sec:brain_tumour_growth}
\citet{mantzavinos_2014} and \citet{asvestas_2014} both consider a one-dimensional reaction-diffusion model for the growth of brain tumours. Due to the heterogeneity of brain tissue, the diffusion coefficient is assumed to be piecewise constant over several regions consisting of either white or grey matter. The spread of malignant cells is governed by the equation \cite{mantzavinos_2014}:
\begin{subequations}
\label{eq:brain_model}
\begin{gather}
\frac{\partial c_{i}}{\partial t} = D_{i}\frac{\partial^{2} c_{i}}{\partial x^{2}} + c_{i},
\end{gather}
for $x\in(l_{i-1},l_{i})$, $t > 0$ and $i = 1,\hdots,m$, where $c_{i}(x,t)$ is the cell density at position $x$ and time $t$ in region $i$ and $D_{i}$ is the constant dimensionless diffusion coefficient in region $i$. Initially, a spatial distribution of malignant cells $f_{i}(x)$ is assumed in each layer: 
\begin{gather}
c_{i}(x,0) = f_{i}(x),
\end{gather}
for $x\in[l_{i-1},l_{i}]$ and $i = 1,\hdots,m$. Continuity of the cell density and flux are assumed at the interfaces between adjacent regions and the migration of cells beyond the brain boundaries is prohibited \cite{mantzavinos_2014}, yielding the following forms for internal and external BCs:
\begin{gather}
 c_{i}(l_{i},t) = c_{i+1}(l_{i},t),\quad 
 D_{i}\frac{\partial c_{i}}{\partial x}(l_{i},t) = D_{i+1}\frac{\partial c_{i+1}}{\partial x}(l_{i},t),\\
\frac{\partial c_{1}}{\partial x}(l_{0},t) = 0,\quad\frac{\partial c_{m}}{\partial x}(l_{m},t) = 0,
\end{gather}
\end{subequations}
for $t > 0$ and $i = 1,\hdots,m-1$. The equation system (\ref{eq:brain_model}) is converted into the required form of the multilayer diffusion problem considered in this paper (i.e., without the source term) via the substitution $c_{i}(x,t) = e^{t}u_{i}(x,t)$ \cite{mantzavinos_2014}, which yields the multilayer diffusion problem described by equations (\ref{eq:original_equation})--(\ref{eq:original_bc2}) and (\ref{eq:ic_implicit}) with $a_{L} = a_{R} = 0$, $b_{L} = b_{R} = 1$ and $g_{0}(t) = g_{m}(t) = 0$.

\begin{figure}[htb]
\centering
\includegraphics[width = 0.55\textwidth]{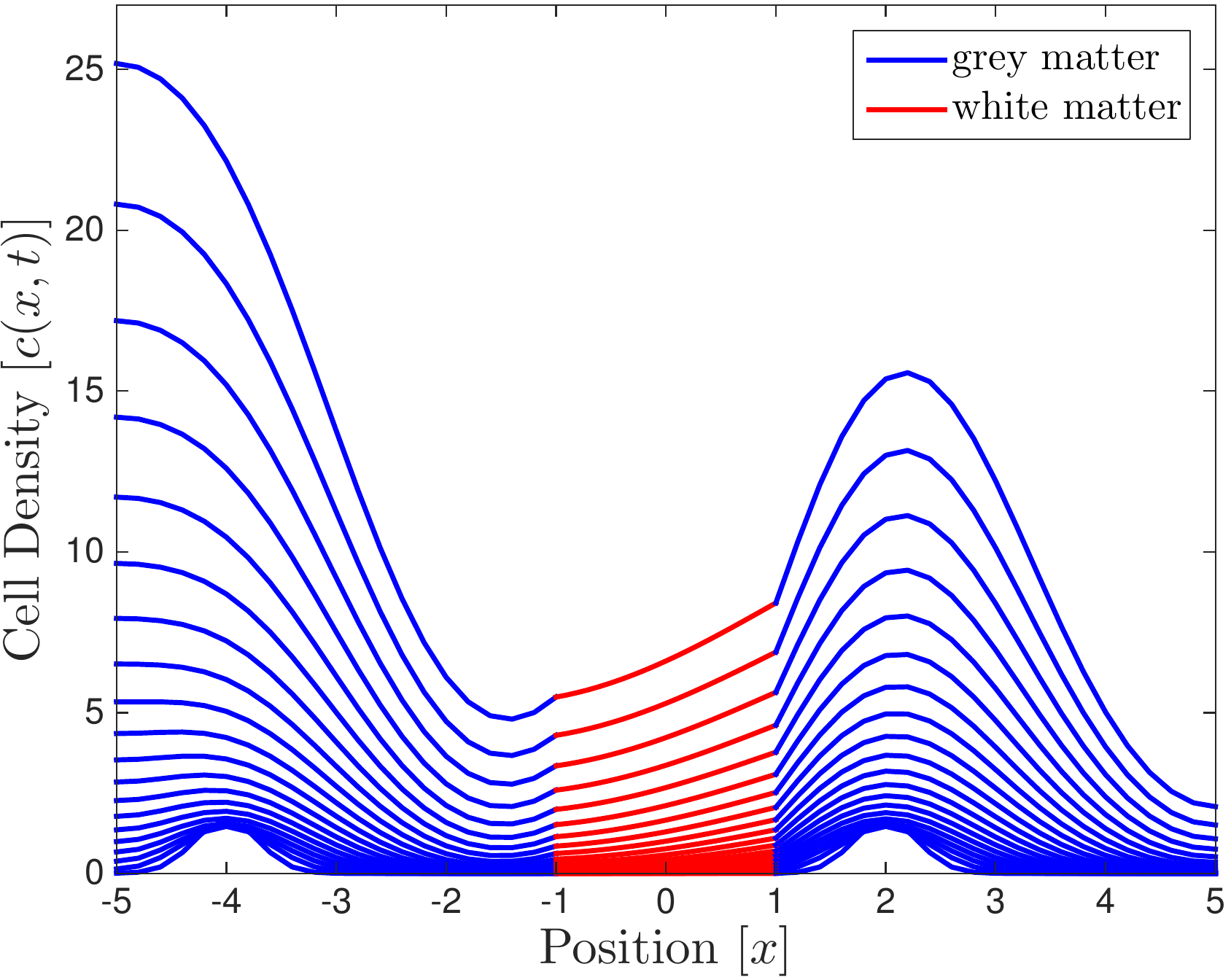}
\caption{Dimensionless cell density $c(x,t)$ distribution for the brain tumour growth problem described in Section \ref{sec:brain_tumour_growth}. Solution is given for $t = 0.2\!:\!0.2\!:\!4.0$.}
\label{fig:mantzavinos}
\end{figure}

We consider the test case described by \citet{mantzavinos_2014} with 
$m = 3$ regions, domain $[l_{0},l_{1},l_{2},l_{3}] = [-5,-1,1,5]$, diffusion coefficients $D_{1} = D_{3} = 0.2$ (grey matter) and $D_{2} = 1$ (white matter), and initial point sources of tumour cells at $x = -4$ and $x = 2$:
\begin{gather*}
u_{1}(x,0) = f_{1}(x) := \delta(x+4),\quad x\in[-5,-1],\\
u_{2}(x,0) = f_{2}(x) := 0,\quad x\in[-1,1],\\
u_{3}(x,0) = f_{3}(x) := \delta(x-2),\quad x\in[1,5],
\end{gather*}
where $\delta(\cdot)$ is the Dirac delta function. To treat the above initial conditions in our code, we use the approximation:
\begin{align*}
\delta(x+\mu) \approx \frac{1}{a\sqrt{\pi}}\exp\left(-\frac{(x-\mu)^{2}}{a^{2}}\right),
\end{align*}
with $a = 0.1$. Figure \ref{fig:mantzavinos} plots the cell density over time computed via the semi-analytical solution, replicating the plot presented in \cite{mantzavinos_2014}.

\section{Conclusions}
This paper has developed a new semi-analytical method for solving  multilayer diffusion problems with time-varying external BCs and general internal BCs at the interfaces between adjacent layers. Numerical experiments suggested the semi-analytical solution exhibits a convergence rate of $O(N^{-3})$, where $N$ is the number of eigenvalues used in the solution expansions, a result which was also confirmed analytically. While this is some way away from the exponential convergence rate of the classical analytical solution, the semi-analytical approach possesses other clear advantages such as requiring only simple eigenvalues and performing well for a large number of layers; both due to it avoiding the solution of a complex transcendental equation for the eigenvalues. Numerical experiments demonstrated that solutions can be obtained using the new semi-analytical method that are more accurate and efficient than Sheils' \cite{sheils_2016} unified transform method while also not exhibiting faulty behaviour at the external boundaries. Finally, in contrast to classical analytical solutions and the unified transform method, only the semi-analytical method introduced in this paper is able to correctly treat problems with both time-dependent external BCs and a large number of layers.

\section*{Acknowledgments}
EJC acknowledges funding from the Australian Research Council (DE150101137).

\appendix
\section{Psi functions}
\label{app:psi_functions}
The functions $\psi_{i,1}(x)$ and $\psi_{i,2}(x)$ described in equation (\ref{eq:wi}) are defined as follows:
\begin{itemize}
\item First Layer ($i=1$)
\begin{align*}
\psi_{1,1}(x) &= 
\begin{cases} \frac{x(x-2l_{1})}{2b_{L}(l_{1}-l_{0})} & \text{for $a_{L}=0$}\\[0.2cm] \frac{1}{a_{L}} & \text{for $a_{L}\neq 0$}\\   
\end{cases}\qquad
\psi_{1,2}(x) = \begin{cases} \frac{x(x-2l_{0})}{2\gamma_{1}(l_{1}-l_{0})} & \text{for $a_{L} = 0$}\\[0.2cm]
\frac{a_{L}(x-l_{0})+b_{L}}{\gamma_{1}a_{L}} & \text{for $a_{L}\neq 0$} \end{cases}
\end{align*}
\item Middle Layers ($i = 2,\hdots,m-1$)
\begin{align*}
\psi_{i,1}(x) &= \frac{x(2l_{i}-x)}{2\gamma_{i}(l_{i}-l_{i-1})}\qquad
\psi_{i,2}(x) = \frac{x(x-2l_{i-1})}{2\gamma_{i}(l_{i}-l_{i-1})}
\end{align*}
\item Last Layer ($i = m$)
\begin{align*}
\psi_{m,1}(x) &= \begin{cases} \frac{x(2l_{m}-x)}{2\gamma_{m}(l_{m}-l_{m-1})} & \text{for $a_{R} = 0$}\\[0.2cm] \frac{a_{R}(x-l_{m})-b_{R}}{\gamma_{m}a_{R}} & \text{for $a_{R}\neq 0$} \end{cases}\qquad
\psi_{m,2}(x) = \begin{cases} \frac{x(x-2l_{m-1})}{2b_{R}(l_{m}-l_{m-1})} & \text{for $a_{R}=0$}\\[0.2cm] \frac{1}{a_{R}} & \text{for $a_{R}\neq 0$}\end{cases}
\end{align*}
\end{itemize}

\section{Linear system}
\label{app:linear_system}
This appendix formulates the entries of the matrix $\mathbf{A}(s)$ and $\mathbf{b}(s)$ featuring in the tridiagonal matrix system (\ref{eq:linear_system}). Let $a_{i,j}$ denote the $(i,j)$ entry of $\mathbf{A}$ and $b_{i}$ denote the $i$th entry of $\mathbf{b}$.

Consider equation (\ref{eq:cins}) and define $c_{i,n}^{(k)}(s)$ ($k = 1,2,3$) as follows:
\begin{multline*}
\overline{c}_{i,n}(s) = \underbracket[0.5pt]{\frac{\beta_{i,5,n}}{s + D_{i}\lambda_{i,n}^{2}}}_{c_{i,n}^{(1)}(s)} + \underbracket[0.5pt]{\left(\frac{D_{i}(\beta_{i,3,n}+\lambda_{i,n}^{2}\beta_{i,1,n})}{s+D_{i}\lambda_{i,n}^{2}}-\beta_{i,1,n}\right)}_{c_{i,n}^{(2)}(s)}\overline{g}_{i-1}(s)\\ + \underbracket[0.5pt]{\left(\frac{D_{i}(\beta_{i,4,n}+\lambda_{i,n}^{2}\beta_{i,2,n})}{s+D_{i}\lambda_{i,n}^{2}}-\beta_{i,2,n}\right)}_{c_{i,n}^{(3)}(s)}\overline{g}_{i}(s).
\end{multline*}
Substituting (\ref{eq:laplace_tranform_ui}) into (\ref{eq:constraint_laplace}), and rearranging identifies
\begin{itemize}
\item the subdiagonal of $\mathbf{A}$:
\begin{gather*}
a_{i,i-1} = \psi_{i,1}(l_{i}) + \sum_{n=0}^{N-1}c_{i,n}^{(2)}(s)\widehat{\phi}_{i,n}(l_{i}),\quad i = 2,\hdots,m-1,
\end{gather*}
\item the diagonal of $\mathbf{A}$:
\begin{multline*}
a_{i,i} = \psi_{i,2}(l_{i}) - \theta_{i}\psi_{i+1,1}(l_{i}) + \frac{1}{H_{i}} + \sum_{n=0}^{N-1}\left[c_{i,n}^{(3)}(s)\widehat{\phi}_{i,n}(l_{i})-\theta_{i}c_{i+1,n}^{(2)}(s)\widehat{\phi}_{i+1,n}(l_{i})\right],\\ i = 1,\hdots,m-1,
\end{multline*}
\item the superdiagonal of $\mathbf{A}$:
\begin{align*}
a_{i,i+1} = -\theta_{i}\psi_{i+1,2}(l_{i})-\theta_{i}\sum_{n=0}^{N-1}c_{i+1,n}^{(3)}(s)\widehat{\phi}_{i+1,n}(l_{i}),\quad i = 1,\hdots,m-2.
\end{align*}
\end{itemize}
The entries of $\mathbf{b}$ are given by:
\begin{align*}
&b_{1} = \sum_{n=0}^{N-1}\left[\theta_{1}c_{2,n}^{(1)}(s)\widehat{\phi}_{2,n}(l_{1}) - c_{1,n}^{(1)}(s)\widehat{\phi}_{1,n}(l_{1})\right] - \left[\psi_{1,1}(l_{1}) + \sum_{n=0}^{N-1}c_{1,n}^{(2)}(s)\widehat{\phi}_{1,n}(l_{1})\right]
\overline{g}_{0}(s),\\
&b_{i} = \sum_{n=0}^{N-1}\left[\theta_{i}c_{i+1,n}^{(1)}(s)\widehat{\phi}_{i+1,n}(l_{i}) - c_{i,n}^{(1)}(s)\widehat{\phi}_{i,n}(l_{i})\right],\quad i =2,\hdots,m-2,\\
&b_{m-1} = \sum_{n=0}^{N-1}\Bigl[\theta_{m-1}c_{m,n}^{(1)}(s)\widehat{\phi}_{m,n}(l_{m-1}) - c_{m-1,n}^{(1)}(s)\widehat{\phi}_{m-1,n}(l_{m-1})\Bigr]\\ &\hspace{4.0cm}+ \theta_{m-1}\left[\psi_{m,2}(l_{m-1})+\sum_{n=0}^{N-1}c_{m,n}^{(3)}(s)\widehat{\phi}_{m,n}(l_{m-1})\right]\overline{g}_{m}(s).
\end{align*}

\section{Analytical solution}
\label{app:analytical_solution}
Using separation of variables, the following analytical solution of the multilayer diffusion problem described by equations (\ref{eq:original_equation})--(\ref{eq:original_bc2}) and (\ref{eq:ic_general}) can be derived \cite{trefry_1999,carr_2016a,hickson_2009a} when $g_{0}$ and $g_{m}$ in equations (\ref{eq:original_bc1})--(\ref{eq:original_bc2}) are independent of time $t$:
\begin{align*}
u_{i}(x,t) = w_{i}(x) + \sum_{n=0}^{\infty}c_{n}e^{-t\lambda_{n}^2}\phi_{i,n}(x),
\end{align*}
where $w_{i}(x)$ denotes the steady-state solution in the $i$th layer (see, e.g., \citet{carr_2016a}).
The eigenvalues ($\lambda$){, which form a global set valid across all layers,} and the non-normalized eigenfunctions ($\phi_{i}$) satisfy a series of coupled Sturm Liouville problems involving homogeneous versions of the internal and external BCs:
\begin{gather*}
-D_{i}\phi_{i}'' = \lambda^{2}\phi_{i},\quad x\in(l_{i-1},l_{i}),\quad i = 1,\hdots, m,\\
a_{L}\phi_{1}(l_{0}) - b_{L}\phi_{1}'(l_{0}) = 0,\\
a_{R}\phi_{m}(l_{m}) + b_{R}\phi_{m}'(l_{m}) = 0,\\
\gamma_{i}\phi_{i}'(l_{i}) = H_{i}(\theta_{i}\phi_{i+1}(l_{i}) - \phi_{i}(l_{i})),\quad i=1,\hdots,m-1,\\
\gamma_{i}\phi_{i}'(l_{i}) = \gamma_{i+1}\phi_{i+1}'(l_{i}),\quad i=1,\hdots,m-1.
\end{gather*}
Substituting the form of the eigenfunctions:
\begin{gather*}
\phi_{i}(\lambda; x) = \zeta_{i}(\lambda)\sin\left(\frac{\lambda}{\sqrt{D_{i}}}(x-l_{i-1})\right) + \xi_{i}(\lambda)\cos\left(\frac{\lambda}{\sqrt{D_{i}}}(x-l_{i-1})\right),
\end{gather*}
into the above internal and external BCs yields a linear system, which can be expressed in matrix form as $\mathbf{A}(\lambda)\mathbf{x} = \mathbf{0}$, where $\mathbf{x} = (\zeta_{1},\xi_{1},\hdots,\zeta_{m},\xi_{m})^{T}$ and $\mathbf{A}\in\mathbb{R}^{2m\times 2m}$. The notation $\zeta_{i}(\lambda)$, $\xi_{i}(\lambda)$ and $\mathbf{A}(\lambda)$ is used since they each depend on $\lambda$. The eigenvalues $\lambda_{n}$ ($n = 0,1,\hdots$) are the non-negative roots of the transcendental equation:
\begin{align*} 
\det(\mathbf{A}(\lambda)) = 0.
\end{align*}
For each eigenvalue $\lambda_{n}$ ($n = 0,1,\hdots$), an eigenfunction $\phi_{i,n}(x) := \phi_{i}(\lambda_{n}; x)$ is defined, where the coefficients $\zeta_{i}(\lambda_{n})$ and $\xi_{i}(\lambda_{n})$ ($i = 1,\hdots,m$) are determined by finding a non-trivial solution of $\mathbf{A}(\lambda_{n})\mathbf{x} = \mathbf{0}$ \cite{carr_2016a}. The eigenfunctions are orthogonal over the full domain $[l_{0},l_{m}]$ with respect to the weight function $p_{i}(x) = \gamma_{i}D_{i}^{-1}\prod_{k=1}^{i-1}\theta_{k}$. A proof of this result follows closely the proof given by \citet{trefry_1999} for $\gamma_{i}=D_{i}$. Hence, using the initial condition, we have that:
\begin{align*}
c_{n} = \frac{\sum_{i=1}^{m}\int_{l_{i-1}}^{l_{i}} p_{i}(x)\tilde{f}_{i}(x)\phi_{i,n}(x)\,dx}{\sum_{i=1}^{m}\int_{l_{i-1}}^{l_{i}}p_{i}(x)\phi_{i,n}^{2}(x)\,dx},
\end{align*}
where $\tilde{f}_{i}(x) = f_{i}(x)-w_{i}(x)$.

\bibliographystyle{plainnat}
\bibliography{references}

\end{document}